\RequirePackage{ifpdf}
\ifpdf 
\documentclass[pdftex]{sigma}
\else
\documentclass{sigma}
\fi

\newcommand{\x}{{\boldsymbol x}}
\newcommand{\y}{{\boldsymbol y}}
\newcommand{\gf}{{\boldsymbol g}}
\newcommand{\h}{{\boldsymbol h}}
\newcommand{\uu}{{\boldsymbol u}}
\newcommand{\vv}{{\boldsymbol v}}
\newcommand{\w}{{\boldsymbol w}}
\newcommand{\g}{{\mathcal G}}
\newcommand{\z}{{\mathcal Z}}

\begin{document}
\allowdisplaybreaks

\renewcommand{\PaperNumber}{039}

\FirstPageHeading

\ShortArticleName{Combined Reduced-Rank Transform}

\ArticleName{Combined Reduced-Rank Transform}

\Author{Anatoli  TOROKHTI and Phil HOWLETT}
\AuthorNameForHeading{A. Torokhti and P. Howlett}

\Address{School of Mathematics and Statistics, University of South Australia,
Australia}
\Email{\href{mailto:anatoli.torokhti@unisa.edu.au}{anatoli.torokhti@unisa.edu.au},
\href{mailto:phil.howlett@unisa.edu.au}{phil.howlett@unisa.edu.au}}
\URLaddress{\url{http://people.unisa.edu.au/Anatoli.Torokhti}\\
\hspace*{10.5mm}\url{http://people.unisa.edu.au/Phil.Howlett}}

\ArticleDates{Received November 25, 2005, in f\/inal form March 22,
2006; Published online April 07, 2006}

\Abstract{We propose and justify a new approach to constructing  optimal
nonlinear  transforms of random vectors.
We show that the proposed  transform improves such characteristics of {rank-reduced} transforms
 as  compression ratio,  accuracy of decompression and reduces required computational work.
 The proposed transform ${\mathcal T}_p$ is presented in the form of a sum with $p$ terms  where each
 term is interpreted as  a particular {rank-reduced} transform. Moreover,  terms in ${\mathcal T}_p$
  are represented as a combination of three operations ${\mathcal F}_k$, ${\mathcal Q}_k$ and ${\boldsymbol{\varphi}}_k$ with
  $k=1,\ldots,p$.
 The prime idea is to determine ${\mathcal F}_k$ separately, for each $k=1,\ldots,p$, from an associated
 rank-constrained minimization problem similar to that used  in the  Karhunen--Lo\`{e}ve
 transform.  The operations  ${\mathcal Q}_k$ and ${\boldsymbol{\varphi}}_k$ are auxiliary for f\/inding  ${\mathcal F}_k$. The
 contribution of each term  in  ${\mathcal T}_p$   improves the entire transform performance.
A corresponding {unconstrained}  nonlinear optimal transform is also considered. Such a
transform is important in its own right because it is treated as an optimal f\/ilter  without
signal compression.
A rigorous analysis of errors associated with the proposed transforms is given.}


\Classification{41A29}

\section{Introduction}
\label{intr}

Methods of data dimensionality reduction \cite{hot1,kar1,loe1,jol1,sch1,yam1,hua1,vap1,oca1,
tip1,tip2,scho1,ten1,row1,cri1,yam2,hua2,kne1,hon1,che1,hon2,sto1,fuk1,kra1} have been applied successfully
 to many applied problems. The diversity of  applications has stimulated a considerable
 increase in the study of data dimensionality reduction  in  recent decades.  Signif\/icant
 recent results in this  challenging research area are described, in particular, in references
 \cite{jol1,sch1,yam1,hua1,vap1,oca1,tip1,tip2,scho1,ten1,row1,cri1,yam2,hua2,kne1,hon1,che1,hon2,sto1,fuk1,kra1}.
 The known methods concern  both a probabilistic setting
 (as in \cite{sch1,yam1,hua1,vap1,oca1,tip1,yam2,hua2,kne1,hon1,che1,hon2,sto1,fuk1,kra1,
 tor101,tor1,tor3,tor2}) and  deterministic setting
 (as in \cite{scho1,ten1,row1,cri1}) in the dimensionality reduction. The associated techniques
 are often based on the use of reduced-rank operators.

In this paper, a further advance in the development of reduced-rank transforms  is presented.
We study a new approach to  data dimensionality reduction  in a {\em probabilistic} setting
based on the development of  ideas presented in   \cite{sch1,yam1,hua1,tor1,tor3,tor2,tor4}.

Motivation for the proposed approach arises from the following observation.
In general, the reduced-rank transform consists of the three companion operations which are
f\/iltering, compression and reconstruction \cite{sch1,yam1,hua1,yam2,tor1}.
Filtering and compression are performed simultaneously to estimate a reference signal
$\x$ with $m$ components from noisy observable data $\y$ and to f\/ilter and reduce the data
to a shorter vector $\hat{\x}$ with $\eta$ components, $\eta < m$. Components of  $\hat{\x}$
are often called  principal components \cite{jol1}.  The quotient $\eta/m$ is called the
compression ratio.  Reconstruction returns a vector $\tilde{\x}$ with $m$ components so that
 $\tilde{\x}$ should be close to the original $\x$. It is natural to perform these three
 operations so that the reconstruction error and the related computational burden are minimal.

As a result, the performance of the reduced-rank transform is characterized by  three issues
which are (i)  associated  accuracy, (ii) compression ratio, and (iii) {computational
work}.

For a given compression ratio,  the Karhunen--Lo\`{e}ve transform (KLT) \cite{sch1,yam1,hua1}
minimizes the reconstruction error   over the class of all {\em linear}  reduced-rank
transforms. Nevertheless, it may happen that the accuracy  and  compression ratio associated
with the KLT are still not satisfactory. In such a case, an  improvement in the accuracy and
compression ratio can be achieved by a transform with a more general structure than that of
the KLT. Special {\em non-linear}  transforms have been studied in \cite{yam2,hua2,kne1,hon1,che1,hon2,sto1,fuk1,kra1,
 tor101,tor1,tor3,tor2,tor4,tor44,tor444,son1,che2,how2}
using transform structures  developed from the generalised Volterra polynomials. Nevertheless,
the transforms \cite{yam2,tor1,tor3,tor2,tor4} imply a substantial  computational burden
associated with the large number $N$ of terms required by the underlying Volterra polynomial
structure.

{Our objective} is to justify a new transform that may have both accuracy and compression ratio
better  than those of the known  transforms   \cite{sch1,yam1,hua1,tor1,tor3,tor2,tor4}.
A related objective is to f\/ind a way to reduce the associated   computational work compared
with that implied by the transforms  \cite{tor1,tor3,tor2,tor4}. The  analysis  of these issues
is given in Sections \ref{stat}, \ref{det-fk} (Remark~\ref{remark4}), \ref{sec5.2.3} and \ref{spec}.

In Section \ref{part-cas}, we show that the proposed approach generalizes  the Fourier series
in Hilbert space, the Wiener f\/ilter, the Karhunen--Lo\`{e}ve transform and the transforms given
in   \cite{tor1,tor3,tor4}.

\section{Method description}
\label{summ}

We use the following notation:

 $(\Omega, \Sigma, \mu)$ is a  probability space,
where $\Omega = \{\omega \}$ is the set of outcomes, $\Sigma$ a $\sigma$-f\/ield of
measurable subsets  of $\Omega$ and $\mu:\Sigma \rightarrow [0,1]$ an
associated probability measure on $\Sigma$ with $\mu(\Omega) = 1$;
  ${\x} \in L^{2}(\Omega,{\mathbb R}^{m})$ and  ${\y} \in
L^{2}(\Omega,{\mathbb R}^{n})$
are random vectors with realizations $x=\x(\omega) \in {\mathbb R}^{m}$ and
$y=\y(\omega) \in {\mathbb R}^{n}$, respectively.

Each matrix $M \in
{\mathbb R}^{m \times n}$ def\/ines a bounded linear transformation ${\cal
M}: L^{2}(\Omega,{\mathbb R}^{n})\rightarrow$ $L^{2}(\Omega,{\mathbb R}^{m})$
 via the formula $[{\cal M}{\y}](\omega) = M {\y}(\omega)$
for each $\omega \in \Omega$. We note that there are many bounded linear transformations from
$L^{2}(\Omega,{\mathbb R}^{n})$ into $L^{2}(\Omega,{\mathbb R}^{m})$
that cannot be written in the form $[{\cal M}\y](\omega)
= M\y(\omega)$ for each $\omega \in \Omega$. A trivial example is
${\mathcal A}: L^{2}(\Omega,{\mathbb R}^{n})\rightarrow L^{2}(\Omega,{\mathbb R}^{m})$ given by ${\displaystyle {\mathcal A}(\y) = \int_\Omega \y (\omega) d\mu (\omega). }$

Throughout the paper, the calligraphic character letters denote operators def\/ined similarly
to ${\cal M}$.

Let $\gf=[\gf_1 \ldots \gf_m]^T\in L^2(\Omega, {\mathbb R}^m)$ and
   $\h = [\h_1 \ldots \h_n]^T\in L^2(\Omega, {\mathbb R}^n)$ be
 random vectors with $\gf_i, \h_k\in L^2(\Omega, {\mathbb R})$
for  $i=1,\ldots, m$, $k=1,\ldots,n$.
For all $i=1,\ldots, m$ and  $k=1,\ldots, n$, we set
\begin{gather}
  \label{e1}
E[\gf_i] =  \int_{\Omega} \gf_i(\omega)d\mu (\omega),
 \qquad E[\gf_i \h_k] =  \int_{\Omega} \gf_i(\omega) \h_k(\omega) d\mu (\omega),
 \\
 E_{gh} = E\big[\gf\h^T\big] = \{E[\gf_i \h_k]\}\in {\mathbb R}^{m\times n}\qquad\mbox{and}\qquad
 E_{g} = E[\gf] = \{E[\gf_i]\}\in {\mathbb R}^{m}.
  \end{gather}
We also  write
\[
{\mathbb E}_{gh} = E\big[(\gf - E_g)(\h - E_h)^T\big] =
E_{gh} - E[\gf]E\big[\h^T\big].
\]

Achievement of the above objectives  is based on the presentation of the proposed transform in
the form of a sum with $p$ terms (\ref{t3}) where each term is interpreted as  a particular
rank-reduced transform. Moreover,  terms in (\ref{t3}) are represented as a combination of
three operations~${\mathcal F}_k$, ${\mathcal Q}_k$ and ${\boldsymbol{\varphi}}_k$
for each  $k=1,\ldots,p$, where ${\boldsymbol{\varphi}}_k$ is nonlinear.
 The prime idea is to determine~${\mathcal F}_k$ separately, for each $k=1,\ldots,p$, from an associated
 rank-constrained minimization problem similar to that in the KLT.  The operations  ${\mathcal Q}_k$ and
 ${\boldsymbol{\varphi}}_k$ are auxiliary for f\/inding  ${\mathcal F}_k$.
 It is natural to expect that a contribution of
 each term  in (\ref{t3}) will improve the entire transform performance.

To realize such a scheme, we choose the  ${\mathcal Q}_k$
as orthogonal/orthonormal operators (see Section~\ref{struc}).
Then each  ${\mathcal F}_k$ can be determined independently for each
 individual  problem (\ref{min-k}) or (\ref{min-k2}) below.
Next,  operators  ${\boldsymbol{\varphi}}_k$ are used to
reduce the number of terms from $N$ (as in \cite{yam2,tor1,tor3,tor2,tor4})
to $p$ with $p\ll N$. For example, this can  be done when we choose    ${\boldsymbol{\varphi}}_k$
in the form presented in Section \ref{spec}.
Moreover, the composition of operators  ${\mathcal Q}_k$
and ${\boldsymbol{\varphi}}_k$ allows us to reduce the related covariance
matrices to the identity matrix or to a block-diagonal form  with small blocks.
Remark~\ref{remark4} in Section \ref{det-fk} gives more details in this regard.
The computational work associated with such blocks is much less than
that for the large  covariance matrices in \cite{yam2,tor1,tor3,tor2,tor4}.

To regulate  accuracy associated with the proposed transform and its compression
ratio, we formulate the problem in the form (\ref{min1})--(\ref{con1})
where (\ref{con1}) consists of $p$ constraints. It is shown in Remark~\ref{remark2}
of Section~\ref{stat}, and in Sections \ref{sec5.2.1}, \ref{det-fk} and \ref{spec}
 that such a combination of constraints  allows us to equip
 the proposed transforms with several degrees of freedom.

The structure of our transform is presented in Section \ref{struc} and the formal statement
of the problem in Section~\ref{stat}. In Section  \ref{sol}, we determine operators  ${\mathcal Q}_k$
and ${\mathcal F}_k$ (Lemmata \ref{lemma1} and \ref{ort3}, and  Theorems \ref{sol1} and \ref{sol2}, respectively).

\section{Structure of the proposed transform}
\label{struc}

\subsection{Generic form}

The proposed transform ${\mathcal T}_p$ is presented in the form
\begin{gather}
  \label{t3}
{\mathcal T}_p(\y) = f + \sum _{k=1}^p {\mathcal F}_k{\mathcal Q}_k{\boldsymbol{\varphi}}_k(\y)
= f +{\mathcal F}_1{\mathcal Q}_1{\boldsymbol{\varphi}}_1(\y) + \cdots
+ {\mathcal F}_p{\mathcal Q}_p{\boldsymbol{\varphi}}_p(\y),
\end{gather}
where $f\in {\mathbb R}^m$,  ${\boldsymbol{\varphi}}_k:L^{2}
(\Omega,{\mathbb R}^{n})\rightarrow L^{2}(\Omega, {\mathbb R}^n)$,
${\mathcal Q}_1,\ldots,{\mathcal Q}_p:L^2(\Omega, {\mathbb R}^n)\rightarrow L^2(\Omega, {\mathbb R}^n)$
and  ${\mathcal F}_k: L^{2}(\Omega, {\mathbb R}^n) \rightarrow L^2(\Omega, {\mathbb R}^m)$.

In general, one can put  ${\x} \in L^{2}(\Omega, H_X)$, ${\y} \in
L^{2}(\Omega,H_Y)$,  ${\boldsymbol{\varphi}}_k:L^{2}(\Omega,H_Y)\rightarrow L^{2}(\Omega, H_k)$,
 ${\mathcal Q}_k:L^2(\Omega, H_k)\rightarrow L^2(\Omega, \tilde{H}_k)$
 and  ${\mathcal F}_k:L^2(\Omega, \tilde{H}_k)\rightarrow L^2(\Omega, H_X)$
with $H_X$,  $H_Y$,   $H_k$ and $\tilde{H}_k$ separable  Hilbert spaces, and $k=1,\ldots,p$.

In (\ref{t3}), the vector $f$ and operators ${\mathcal F}_1, \ldots, {\mathcal F}_p$
are determined from the minimization problem (\ref{min1})--(\ref{con1})
given in the  Section \ref{stat}.
Operators ${\mathcal Q}_1,\ldots,{\mathcal Q}_p$ in (\ref{t3})
are orthogonal (orthonormal) in the sense of the  Def\/inition~\ref{def1}
in Section \ref{stat} (in this regard, see also Remark \ref{remark3} in Section~\ref{ort}).

To demonstrate and justify f\/lexibility of the transform ${\mathcal T}_p$
with respect to the choice
of ${\boldsymbol{\varphi}}_1,\ldots,{\boldsymbol{\varphi}}_p$
in (\ref{t3}), we mainly study the case   where  ${\boldsymbol{\varphi}}_1,\ldots,{\boldsymbol{\varphi}}_p$
are arbitrary. Specif\/ications of  ${\boldsymbol{\varphi}}_1,\ldots,{\boldsymbol{\varphi}}_p$
are presented in Sections \ref{some}, \ref{spec} and \ref{part-cas} where we also discuss
the benef\/its associated with some particular forms of ${\boldsymbol{\varphi}}_1,\ldots,{\boldsymbol{\varphi}}_p$.

\subsection{Some particular cases}
\label{some}

 Particular cases of the model ${\mathcal T}_p$ are associated
 with specif\/ic choices of  ${\boldsymbol{\varphi}}_k$,  ${\mathcal Q}_k$ and ${\mathcal F}_{ k}$.
 Some examples are given below.

 {\bf (i)} If $H_X = H_Y = {\mathbb R}^{n}$ and
$H_k = \tilde{H}_k = {\mathbb R}^{nk}$ where ${\mathbb R}^{nk}$ is
the $k$th degree of  ${\mathbb R}^{n}$, then (\ref{t3}) generalises
the known transform structures  \cite{yam2,tor1,tor3,tor2,tor4}.
The models  \cite{yam2,tor1,tor3,tor2,tor4} follow from  (\ref{t3})
if   ${\boldsymbol{\varphi}}_k(\y) = \y^k$ where  $\y^k = (\y,\ldots, \y)\in L^{2}(\Omega,{\mathbb R}^{nk})$,
${\mathcal Q}_k = {\cal I}$, where ${\cal I}$ is the identity operator,   and if ${\mathcal F}_k$
is a $k$-linear operator. It has been shown in \cite{yam2,tor1,tor3,tor2,tor4}
that such a form of ${\boldsymbol{\varphi}}_k$ leads to a signif\/icant improvement in the associated accuracy.
 See Section~\ref{part-cas} for more details.

{\bf (ii)} If  ${\boldsymbol{\varphi}}_k:L^{2}(\Omega,H_Y)\rightarrow L^{2}(\Omega, H_X)$
and $\{\uu_1,  \uu_2, \ldots \}$ is a basis in $L^{2}(\Omega, H_X)$
 then ${\boldsymbol{\varphi}}_k$ and~${\mathcal Q}_k$ can be chosen
 so that ${\boldsymbol{\varphi}}_k(\y) = \uu_k$ and  ${\mathcal Q}_k = {\cal I}$, respectively.
As a result, in this particular case, $ {\mathcal T}_p(\y)
= f + \sum\limits_{k=1}^p {\mathcal F}_k(\uu_k).$

{\bf (iii)} A similar case follows if  ${\boldsymbol{\varphi}}_k:L^{2}(\Omega, H_Y)\rightarrow L^{2}(\Omega, H_k)$
is arbitrary but  ${\mathcal Q}_k:L^2(\Omega, H_k)\rightarrow  L^2(\Omega, \tilde{H}_k)$
is def\/ined so that ${\mathcal Q}_k[{\boldsymbol{\varphi}}_k(\y)] = \vv_k$   with $k=1,\ldots , p$
where
 $\{\vv_1, \vv_2, \ldots\}$ is a basis in $L^2(\Omega, \tilde{H}_k)$.
 Then  ${\mathcal T}_p(\y) = f + \sum\limits _{k=1}^p {\mathcal F}_k(\vv_k).$

{\bf (iv)} Let $\tilde{\x}^{(1)},\ldots,\tilde{\x}^{(p)}$ be estimates of $\x$
by the known transforms
\cite{hua1,tor101,tor44}. Then we can put ${\boldsymbol{\varphi}}_1(\y) = \tilde{\x}^{(1)},$
$\ldots,$ ${\boldsymbol{\varphi}}_p(\y) = \tilde{\x}^{(p)}.$ In particular,
one could choose ${\boldsymbol{\varphi}}_1(\y) = \y$. In such a~way, the
vector $\x$ is pre-estimated from $\y$, and therefore, the overall $\x$
estimate by ${\mathcal T}_p$ will be improved.
 A new recursive method for f\/inding
 $\tilde{\x}^{(1)},\ldots,\tilde{\x}^{(p)}$ is given in Section \ref{spec} below.

Other  particular cases of the proposed transform are  considered in Sections \ref{spec} and \ref{part-cas}.

\begin{remark}\label{remark1}
The particular case of ${\mathcal T}_p$ considered in the item {\bf (iii)}
above can be interpreted as an operator form of the Fourier polynomial
in Hilbert space  \cite{cot1}.  The benef\/its associated with the Fourier polynomials
are well known. In item {\bf (ii)} of Section \ref{part-cas}, this case is considered in more detail.
\end{remark}

\section{Statement of the problem}
\label{stat}

First, we def\/ine  orthogonal and orthonormal operators as follows.

\begin{definition}\label{def1}
Let $\uu_k\in L^{2}(\Omega,{\mathbb R}^{n})$ and $\vv_k = {\mathcal Q}_k(\uu_k)$.
The operators ${\mathcal Q}_1,\ldots,{\mathcal Q}_p$ are called  pairwise orthonormal if
$
{\mathbb E}_{v_i v_j} = \left \{ \begin{array}{@{}cc}
{\mathbb O}, &  i\neq j,\\
 I, &  i=j \end{array} \right.
$
 for any  $i,j = 1,\ldots, p$. Here, ${\mathbb O}$ and $I$ are the zero matrix and identity matrix, respectively.
If
$
{\mathbb E}_{v_i v_j} ={\mathbb O}\quad\mbox{for}\quad i\neq j
$
 with  $i,j = 1,\ldots, p$, and if ${\mathbb E}_{v_i v_j}$
 is not necessarily equal to $I$ for $i=j $  then ${\mathcal Q}_1,\ldots,{\mathcal Q}_p$
 are called  pairwise orthogonal.
 \end{definition}

Hereinafter, we suppose that ${\mathcal F}_k$ is linear for all $k=1,\ldots,p$ and that the Hilbert spaces are the f\/inite dimensional Eucledian spaces, $H_X = {\mathbb R}^m$ and  $H_Y =H_k = \tilde{H}_k = {\mathbb R}^n.$
 For any vector
 $\gf\in L^2(\Omega, {\mathbb R}^m)$, we set
\begin{gather}
  \label{e11}
 E[\|\gf\|^2] = \int_{\Omega}\|\gf(\omega)\|^2 d\mu (\omega) < \infty,
 \end{gather}
 where  $\|\gf(\omega)\|$ is the Euclidean norm of $\gf(\omega)$.

Let us denote
\begin{gather}
  \label{j1}
{J}(f,{\mathcal F}_1,\ldots {\mathcal F}_p) = E\big[\|\x - {\mathcal T}_p(\y)\|^2\big].
\end{gather}

The problem is

(i) to f\/ind   operators ${\mathcal Q}_1,\ldots,{\mathcal Q}_p$ satisfying Def\/inition \ref{def1}, and

(ii) to determine the vector $f^0$ and operators
 ${\mathcal F}_1^0, \ldots, {\mathcal F}_p^0$ such that
\begin{gather}
  \label{min1}
{J}(f^0,{\mathcal F}_1^0,\ldots {\mathcal F}_p^0)
=\min_{f, {\cal F}_1,\ldots, {\cal F}_p} {J}(f,{\mathcal F}_1,\ldots, {\mathcal F}_p)
\end{gather}
{\samepage subject to
\begin{gather}
  \label{con1}
\mbox{rank}\,  {\mathcal F}_1 =  \eta_1, \quad \ldots, \quad \mbox{rank}\, {\mathcal F}_p =  \eta_p,
\end{gather}
where $\eta_1 + \cdots + \eta_p = \eta \leq \min \{m,n\}$.}

Here, for $k=1,\ldots,p$, (see, for example,  \cite{kow1})
\[
\mbox{rank} ({\mathcal F}_k) = \dim {\mathcal F}_k (L^2(\Omega,{\mathbb R}^n)).
\]

We write
\begin{gather}
\label{th1}
\displaystyle{{\mathcal T}_p^0(\y)} = \displaystyle{f^0 + \sum _{k=1}^p {\mathcal F}_k^0(\vv_k)}
\end{gather}
with $\vv_k$ def\/ined by Def\/inition \ref{def1}.

It is supposed that covariance matrices formed from
vectors ${\mathcal Q}_1{\boldsymbol{\varphi}}_1(\y), \ldots, {\mathcal Q}_p{\boldsymbol{\varphi}}_p(\y)$
in (\ref{t3})  are known or can be estimated. Various estimation
methods can be found in \cite{per1,kau1,sch3,kub1,led1,leu1}.
 We note that such an assumption is traditional \cite{hot1,kar1,loe1,jol1,sch1,yam1,hua1,vap1,oca1,
tip1,tip2,scho1,ten1,row1,cri1,yam2}
 in the study of optimal transforms. The ef\/fective estimate of covariance
  matrices represents a specif\/ic task \cite{per1,kau1,sch3,kub1,led1,leu1}
   which is not considered in this paper.

\begin{remark}\label{remark2}
Unlike  known rank-constrained problems, we consider $p$ constraints (\ref{con1}).
The number $p$ of the  constraints and the ranks $\eta_1,\ldots,\eta_p$ form
the degrees of freedom for ${\mathcal T}_p^0$. Variation of $p$ and  $\eta_1,\ldots,\eta_p$
allows us to regulate accuracy associated with the transform ${\mathcal T}_p^0$ (see (\ref{er1})
in Section \ref{sec5.2.1} and (\ref{er12}) in Section \ref{det-fk}) and its compression ratio
(see (\ref{cr}) in Section \ref{spec}). It follows from (\ref{er1}) and (\ref{er12})
 that the accuracy increases if $p$ and  $\eta_1,\ldots,\eta_p$ increase.
 Conversely, by (\ref{cr}), the compression ratio is improved if  $\eta_1,\ldots,\eta_p$  decrease.
\end{remark}

\section{Solution of the problem}
\label{sol}

The problem (\ref{min1})--(\ref{con1}) generalises the
known rank-constrained problems  where  only  one constraint has been considered.
Our plan for the solution is as follows. First, in Section \ref{ort},
we will determine the operators ${\mathcal Q}_1,\ldots,{\mathcal Q}_p$.
Then, in Section~\ref{f0fk}, we will obtain  $f^0$ and
 ${\mathcal F}_1^0, \ldots, {\mathcal F}_p^0$
 satisfying (\ref{min1}) and (\ref{con1}).

\subsection[Determination of orthogonalizing operators ${\mathcal Q}_1,\ldots,{\mathcal Q}_p$]{Determination
of orthogonalizing operators $\boldsymbol{{\mathcal Q}_1,\ldots,{\mathcal Q}_p}$}
\label{ort}

If $M$ is a square matrix then we write $M^{1/2}$ for a matrix such that
$
M^{1/2}M^{1/2} = M.
$
 We note that the matrix $M^{1/2}$ can be computed in various ways \cite{hig1}.
 In this paper,   $M^{1/2}$ is determined from the singular value decomposition (SVD) \cite{gol1} of $M$.

For the case when matrix ${\mathbb E}_{v_k v_k}$
is invertible for any $k=1,\ldots,p$, the orthonormalization
procedure is as follows. For $\uu_k\in L^{2}(\Omega,{\mathbb R}^{n})$, we write
\begin{gather}
\label{qk}
[{\mathcal Q}_k(\uu_k)](\omega) = Q_k\uu_k(\omega),
\end{gather}
where $Q_k\in {\mathbb R}^{n\times n}$.
For $\uu_k, \vv_j, \w_j\in L^2(\Omega, {\mathbb R}^n)$,
we also def\/ine operators  ${\mathcal E}_{u_k v_j}, {\mathcal E}^{-1}_{v_j v_j}:
L^{2}(\Omega,{\mathbb R}^{n})\!\rightarrow$ $ L^{2}(\Omega,{\mathbb R}^{n})$ by the equations
\begin{gather}
\label{ge}
[{\mathcal E}_{u_k v_j}(\w_j)](\omega)
= {\mathbb E}_{u_k v_j}\w_j(\omega)\qquad\mbox{and}\qquad
[{\mathcal E}^{-1}_{v_j v_j}(\w_j)](\omega) = {\mathbb E}^{-1}_{v_j v_j}\w_j(\omega),
\end{gather}
respectively.

\begin{lemma}\label{lemma1}  Let
\begin{gather}
\label{w}
\w_1 =\uu_1\qquad\mbox{and}\qquad \displaystyle{\w_i = \uu_i - \sum_{k=1}^{i-1}
{\mathcal E}_{u_i w_k}{\mathcal E}^{-1}_{w_k w_k} (\w_k })\quad\mbox{for}\quad i=1,\ldots,p,
\end{gather}
 where ${\mathcal E}^{-1}_{w_k w_k}$ exists.
Then

(i) the vectors $\w_1,\ldots,\w_p$ are pairwise orthogonal, and

(ii) the vectors $\vv_1,\ldots,\vv_p$, defined by
\begin{gather}
\label{wq}
 \vv_i ={\mathcal Q}_i(\uu_i)
\end{gather}
with
\begin{gather}
\label{wq1}
{\mathcal Q}_i(\uu_i) = \big({\mathcal E}^{1/2}_{w_i w_i}\big)^{-1} (\w_i)
\end{gather}
for $i=1,\ldots,p$,  are pairwise orthonormal.
\end{lemma}

\begin{proof}
The proof is given in the Appendix.
\end{proof}

For the case when matrix ${\mathbb E}_{v_k v_k}$ is singular for $k=1,\ldots,p$,
the orthogonalizing operators  ${{\mathcal Q}}_1, \ldots ,{{\mathcal Q}}_p$
 are determined by  Lemma~\ref{ort3} below.
 Another dif\/ference from Lemma~\ref{lemma1} is that the vectors  $\vv_1,\ldots,\vv_p$
 in Lemma~\ref{ort3}  are pairwise orthogonal but not orthonormal.
An intermediate result is given in Lemma~\ref{phil}.

The symbol $\dag$ is used to denote the pseudo-inverse operator~\cite{ben1}.
It is supposed that the pseudo-inverse $M^\dag$ for matrix $M$ is determined from the SVD of $M$.

\begin{lemma}[\cite{tor1}]
\label{phil}
For any random vectors ${\boldsymbol g}\in L^2(\Omega,
{\mathbb R}^m)$ and  ${\boldsymbol h}\in L^2(\Omega, {\mathbb R}^n)$,
\begin{gather}
\label{egh}
{\mathbb E}_{gh}{\mathbb E}^\dag_{hh}{\mathbb E}_{hh} = {\mathbb E}_{gh}.
\end{gather}
\end{lemma}

\begin{lemma}
\label{ort3}
Let $\vv_i = {\mathcal Q}_i(\uu_i)$ for $i=1,\ldots,p$, where ${{\mathcal Q}}_1, \ldots ,{{\mathcal Q}}_p$ are such that
\begin{gather}
\label{f1}
{{\mathcal Q}}_1(\uu_1) = \uu_1\qquad\mbox{and}\qquad
{{\mathcal Q}}_i(\uu_i) = \uu_i - \sum_{k=1}^{i-1}\z_{ik}(\vv_k)\quad\mbox{for}\quad i=2,\ldots, p
\end{gather}
with $\z_{ik}:L^{2}(\Omega,{\mathbb R}^{n})\rightarrow
L^{2}(\Omega,{\mathbb R}^{n})$  defined by
\begin{gather}
\label{a}
Z_{ik} = {\mathbb E}_{u_i v_k} {\mathbb E}^\dag_{v_k v_k}
+ A_{ik}(I - {\mathbb E}_{v_k v_k}{\mathbb E}^\dag_{v_k v_k})
\end{gather}
with $A_{ik}\in {\mathbb R}^{n\times n}$ arbitrary.
Then  the vectors $\vv_1,\ldots,\vv_p$  are pairwise orthogonal.
\end{lemma}

\begin{proof} The proof is given in the Appendix.
\end{proof}

We note that Lemma~\ref{ort3} does not require invertibility of matrix ${\mathbb E}_{v_k v_k}$.
At the same time, if ${\mathbb E}^{-1}_{v_k v_k}$ exists,
then  vectors $\w_1,\ldots,\w_p$ and  $\vv_1,\ldots,\vv_p$
def\/ined by (\ref{w}) and  Lemma~\ref{ort3} respectively, coincide.

\begin{remark}\label{remark3}
Orthogonalization of random vectors is not, of course,  a new idea.
In particular,  generalizations of the Gram--Schmidt orthogonalization
procedure  have been considered in \cite{mat1,gol2}.
The proposed orthogonalization procedures in Lemmata \ref{lemma1}
and \ref{ort3} are  dif\/ferent from those in \cite{mat1,gol2}. In particular,
Lemma \ref{ort3} establishes the vector orthogonalization  in terms of  pseudo-inverse operators.
A particular case of the practical implementation of the random
vector orthogonalization is considered in Section~\ref{num}.
\end{remark}

\subsection[Determination of $f^0$, ${\mathcal F}_1^0, \ldots, {\mathcal F}_p^0$
satisfying (\ref{min1})--(\ref{con1})]{Determination of $\boldsymbol{f^0}$, $\boldsymbol{{\mathcal F}_1^0,
\ldots, {\mathcal F}_p^0}$ satisfying (\ref{min1})--(\ref{con1})}
\label{f0fk}

\subsubsection[The case when matrix ${\mathbb E}_{v_iv_i}$ is invertible for $i=1,\ldots,p$]{The case
when matrix $\boldsymbol{{\mathbb E}_{v_iv_i}}$ is invertible for $\boldsymbol{i=1,\ldots,p}$}\label{sec5.2.1}

We consider the simpler case when  ${\mathbb E}_{v_iv_i}$ is invertible for all $i=1,\ldots,p$. Then
the vector $f^0$ and operators ${\mathcal F}_1^0, \ldots,
{\mathcal F}_p^0$ satisfying (\ref{min1})--(\ref{con1}) are def\/ined from
the following Theorem~\ref{sol1}. For each $i=1,\ldots,p$,
let $U_i \Sigma_i V^{T}_i$ be the  SVD of ${\mathbb E}_{xv_i}$,
\begin{gather}
\label{svd}
 U_i \Sigma_i V^{T}_i = {\mathbb E}_{xv_i},
\end{gather}
 where
$U_i\in {\mathbb R}^{m\times n}$, $V_i\in {\mathbb R}^{n\times n}$ are  orthogonal
and  $\Sigma_i\in {\mathbb R}^{n\times n}$ is diagonal,
\begin{gather}
\label{sqd1}
U_i = [s_{i1},\ldots,s_{in}],\qquad  V_i =[d_{i1},\ldots,d_{in}]
\qquad\mbox{and}\qquad \Sigma_i = \mbox{diag}\,(\alpha_{i1},\ldots,\alpha_{in})
\end{gather}
 with $\alpha_{i1} \geq
\cdots \geq \alpha_{ir} > 0$, $\alpha_{i,r+1} = \cdots =
\alpha_{in} = 0$   and $r=1,\ldots,n$ where $r=r(i)$.
We set
\[
U_{i \eta_i} = [s_{i1},\ldots,s_{i\eta_i}],\qquad V_{i \eta_i} = [d_{i1},\ldots, d_{i\eta_i}]\qquad
\mbox{and}\qquad \Sigma_{i \eta_i} = \mbox{diag}(\alpha_{i1},\ldots,\alpha_{i\eta_i}),
\]
where $U_{i\eta_i}\in {\mathbb R}^{m\times \eta_i}$, $V_{i\eta_i}\in {\mathbb R}^{n\times \eta_i}$  and  $\Sigma_{i\eta_i}\in {\mathbb R}^{\eta_i \times \eta_i}$.
 Now we  def\/ine $K_{i\eta_i} \in{\mathbb R}^{m\times n}$ and ${\mathcal K}_{i\eta_i}:L^{2}(\Omega,{\mathbb R}^{n})\rightarrow L^{2}(\Omega,{\mathbb R}^{m})$ by
\begin{gather}
\label{trsvd}
K_{i\eta_i}  = U_{i\eta_i}\Sigma_{i \eta_i}V_{i\eta_i}^{T}\qquad\mbox{and}\qquad [{\mathcal K}_{i\eta_i}(\w_i)](\omega) = K_{i\eta_i}[\w_i(\omega)],
\end{gather}
respectively, for any $\w_i\in L^{2}(\Omega,{\mathbb R}^{n})$.

\begin{theorem}
\label{sol1}
Let $\vv_1,\ldots,\vv_p$ be determined by Lemma~{\rm \ref{lemma1}}.
Then the vector $f^0$ and operators ${\mathcal F}_1^0, \ldots, {\mathcal F}_p^0$,
satisfying \eqref{min1}--\eqref{con1}, are determined by
\begin{gather}
\label{sol-f0}
f^0 = E[\x] - \sum_{k=1}^p F^0_k E[\vv_k]\qquad \mbox{and}\qquad
{\mathcal F}_1^0 = {\mathcal K}_{1\eta_1}, \quad \ldots, \quad {\mathcal F}_p^0  =  {\mathcal K}_{p\eta_p}.
\end{gather}
The accuracy associated with transform ${\mathcal T}_p^0$, determined by  \eqref{th1}
 and \eqref{sol-f0}, is given by
\begin{gather}
\label{er1}
E[\|\x - {\mathcal T}_p^0(\y)\|^2] =\|{\mathbb E}_{xx}^{1/2}\|^2 -\sum_{k=1}^p \sum_{j=1}^{\eta_k} \alpha^2_{kj} .
\end{gather}
\end{theorem}

\begin{proof} The functional $J(f,{\mathcal F}_1, \ldots ,{\mathcal F}_p)$ is written as
\begin{gather}
J(f,{\mathcal F}_1, \ldots ,{\mathcal F}_p)  =  \mbox{tr}\Bigg[E_{xx}
- E[\x]f^T - \sum_{i=1}^p E_{xv_i}F_i^T - f E[\x^T] +  ff^T
+   f\sum_{i=1}^p E[\vv_i^T]F_i^T\nonumber \\
\phantom{J(f,{\mathcal F}_1, \ldots ,{\mathcal F}_p)  =}{}
-   \sum_{i=1}^p F_i E_{v_ix}
+ \sum_{i=1}^p F_iE[\vv_i]f^T + E\Bigg(\sum_{i=1}^p {\mathcal F}_i
(\vv_i) \Bigg[\sum_{k=1}^p {\mathcal F}_i (\vv_i)\Bigg]^T\Bigg)\Bigg].  \label{jqa}
\end{gather}
We remind (see Section \ref{summ}) that here and below, $F_i$
is def\/ined by $[{\mathcal F}_i(\vv_i)](\omega) = F_i[\vv_i(\omega)]$
so that, for example, $E[{\mathcal F}_k(\vv_k)\x_k^T] = F_k E_{v_k x_k}$.
In other words, the right hand side in (\ref{jqa}) is a function of $f$, ${\mathcal F}_1, \ldots ,{\mathcal F}_p$.

Let us show that  $J(f,{\mathcal F}_1, \ldots ,{\mathcal F}_p)$ can be represented as
\begin{gather}
\label{j012}
J(f,{\mathcal F}_1, \ldots ,{\mathcal F}_p) = J_0 + J_1 + J_2,
\end{gather}
where
\begin{gather}
\label{j0}
J_0 = \|{\mathbb E}_{xx}^{1/2}\|^2 -\sum_{i=1}^p  \|{\mathbb E}_{xv_i}\|^2,
\\
\label{j12}
J_1 = \|f - E[\x] + \sum_{i=1}^p F_i E[\vv_i]\|^2\qquad\mbox{and}\qquad
J_2 =\sum_{i=1}^p \|F_i - {\mathbb E}_{xv_i}\|^2.
\end{gather}
Indeed, $J_1$ and  $J_2$ are rewritten as follows
\begin{gather}
J_1 = \mbox{tr} \Bigg(ff^T - fE[\x^T] + \sum_{i=1}^p f E[\vv_i^T]F_i + E[\x] E[\x^T]
 -  E[\x]f^T - \sum_{i=1}^p E[\x] E[\vv_i^T] F_i^T \nonumber\\
\phantom{J_1 =}{} +  \sum_{i=1}^p F_i E[\vv_i]f^T
- \sum_{i=1}^p F_i E[\vv_i] E[\x^T] + \sum_{i=1}^p F_i E[\vv_i]  \sum_{k=1}^p E[\vv_k^T] F_k^T\Bigg)\label{j11}
\end{gather}
and
\begin{gather}
\label{j2}
J_2  =  \sum_{i=1}^p \mbox{tr}\, (F_i - {\mathbb E}_{xv_i})(F_i^T - {\mathbb E}_{v_i x})
= \sum_{i=1}^p \mbox{tr} \,(F_i  F_i^T - F_i {\mathbb E}_{v_i x} - {\mathbb E}_{xv_i}F_i^T + {\mathbb E}_{xv_i} {\mathbb E}_{v_ix}).
\end{gather}
In (\ref{j2}), $\sum\limits_{i=1}^p \mbox{tr}\, (F_i  F_i^T)$  can be represented in the form
\begin{gather}
\label{aea}
\sum_{i=1}^p \mbox{tr} \,(F_i  F_i^T) = \mbox{tr}\Bigg[E\Bigg(\sum_{i=1}^p F_i \vv_i \sum_{k=1}^p \vv_k^T F_k^T\Bigg)
\Bigg] - \mbox{tr}\Bigg(\sum_{i=1}^p F_i E[\vv_i]  \sum_{k=1}^p E[\vv_k^T] F_k^T\Bigg)
\end{gather}
because
\begin{gather}
\label{v-ort}
E[\vv_i \vv_k^T] - E[\vv_i] E[ \vv_k^T] = \left \{ \begin{array}{@{}cc}
{\mathbb O}, &  i\neq k,\\
 I, &  i=k \end{array} \right.
\end{gather}
 due to the orthonormality of  vectors $\vv_1, \ldots, \vv_p$.

Then
\begin{gather}
J_0 + J_1 + J_2  =  \mbox{tr}(E_{xx} - E[\x]E[\x^T]) - \sum_{i=1}^p \mbox{tr} [{\mathbb E}_{xv_i}{\mathbb E}_{v_ix}]\label{jj012}\\
\phantom{J_0 + J_1 + J_2  =}{} + \mbox{tr} \Bigg(ff^T - fE[\x^T] + \sum_{i=1}^p f E[\vv_i^T]F_i
+ E[\x] E[\x^T] -  E[\x]f^T  \nonumber\\
\phantom{J_0 + J_1 + J_2  =}{} - \sum_{i=1}^p E[\x] E[\vv_i^T] F_i^T +  \sum_{i=1}^p F_i E[\vv_i]f^T
 - \sum_{i=1}^p F_i E[\vv_i] E[\x^T] \nonumber\\
\phantom{J_0 + J_1 + J_2  =}{} + \sum_{i=1}^p F_i E[\vv_i]  \sum_{k=1}^p E[\vv_k^T] F_k^T\Bigg)
 +  \mbox{tr}\Bigg[E\Bigg(\sum_{i=1}^p F_i \vv_i \sum_{k=1}^p \vv_k^T F_k^T\Bigg)\Bigg]  \nonumber\\
\phantom{J_0 + J_1 + J_2  =}{} -  \mbox{tr}\Bigg(\sum_{i=1}^p F_i E[\vv_i]  \sum_{k=1}^p E[\vv_k^T] F_k^T\Bigg)\nonumber\\
\phantom{J_0 + J_1 + J_2  =}{} -  \sum_{i=1}^p \mbox{tr} (F_i E_{v_i x} - F_i E[\vv_i] E[\x^T] + E_{xv_i}F_i^T-E[\x]E[\vv_i^T]F_i^T - {\mathbb E}_{xv_i} {\mathbb E}_{v_ix})  \nonumber\\
\phantom{J_0 + J_1 + J_2}{} = J(f,{\mathcal F}_1, \ldots ,{\mathcal F}_p).\label{nj012}
\end{gather}
Hence,  (\ref{j012}) is true. Therefore,
\begin{gather}
 J(f,{\mathcal F}_1, \ldots ,{\mathcal F}_p)  =
\|{\mathbb E}_{xx}^{1/2}\|^2 - \sum_{k=1}^p \|{\mathbb E}_{xv_k}\|^2
+ \|f - E[\x] + \sum_{k=1}^p F_k E[\vv_k]\|^2 \nonumber  \\
\phantom{J(f,{\mathcal F}_1, \ldots ,{\mathcal F}_p)  =}{}+   \sum_{k=1}^p \|F_k - {\mathbb E}_{xv_k}\|^2.\label{proof2}
\end{gather}
It follows from  (\ref{proof2}) that the constrained minimum
(\ref{min1})--(\ref{con1}) is achieved if $f=f^0$ with $f^0$ given by  (\ref{sol-f0}),
and if $F_k^0$ is such that
\begin{gather}
\label{min-k}
J_k(F_k^0) = \min_{F_k} J_k(F_k) \qquad \mbox{subject to} \quad \mbox{rank}\, (F_k) = \eta_k,
\end{gather}
where $J_k(F_k) = \|F_k - {\mathbb E}_{xv_k}\|^2$.
The solution to (\ref{min-k}) is given \cite{gol1} by
\begin{gather}
\label{proof3}
F_k^0 = K_{k\eta_k}.
\end{gather}
Then
\[
E[\|\x - {\mathcal T}_p^0(\y)\|^2]
=\|{\mathbb E}_{xx}^{1/2}\|^2 -\sum_{k=1}^p (\|{\mathbb E}_{xv_k}\|^2  - \|K_{k\eta_k} - {\mathbb E}_{xv_k}\|^2).
\]
Here \cite{gol1},
\begin{gather}
\label{kr}
 \|{\mathbb E}_{xv_k}\|^2 = \sum_{j=1}^{r} \alpha^2_{kj}\qquad\mbox{and}\qquad  \|K_{k\eta_k} - {\mathbb E}_{xv_k}\|^2 =
  \sum_{j=\eta_k +1}^{r} \alpha^2_{kj}
\end{gather}
with $r=r(k)$. Thus, (\ref{er1}) is true. The theorem is proved.
\end{proof}

\begin{corollary}\label{corollary1}
Let $\vv_1,\ldots,\vv_p$ be determined by Lemma~{\rm \ref{lemma1}}.
Then the  vector $\hat{f}$ and operators $\hat{\mathcal F}_1, \ldots, \hat{\mathcal F}_p$
satisfying the unconstrained problem \eqref{min1}, are determined by
\begin{gather}
\label{sol-q}
\hat{f} = E[\x] - \sum_{k=1}^p \hat{F}_k E[\vv_k]\qquad \mbox{and}\qquad
\hat{\mathcal F}_1 = {\mathcal E}_{x v_1}, \quad \ldots, \quad \hat{\mathcal F}_p  =  {\mathcal E}_{x v_p}
\end{gather}
with $\hat{\mathcal F}_k$ such that $[\hat{\mathcal F}_k(\vv_k)](\omega) = \hat{F}_k\vv_k(\omega)$ where
$\hat{F}_k\in{\mathbb R}^{n\times m}$ and $k=1,\ldots,p$.

The accuracy associated with transform $\hat{{\mathcal T}}_p$ given by
\begin{gather}
\label{fr1}
\hat{{\mathcal T}}_p(\y)= \hat{f} + \sum _{k=1}^p \hat{\mathcal F}_k(\vv_k)
\end{gather}
is such that
\begin{gather}
\label{er2}
E[\|\x - \hat{{\mathcal T}}_p(\y)\|^2] =\|{\mathbb E}_{xx}^{1/2}\|^2 -\sum_{k=1}^p \|{\mathbb E}_{xv_k}\|^2 .
\end{gather}
\end{corollary}

\begin{proof} The proof follows directly from (\ref{proof2}).
\end{proof}

\subsubsection[The case when matrix ${\mathbb E}_{v_k v_k}$ is not invertible for $k=1,\ldots,p$]{The case
when matrix $\boldsymbol{{\mathbb E}_{v_k v_k}}$ is not invertible for $\boldsymbol{k=1,\ldots,p}$}
\label{det-fk}

We write $A_k\in {\mathbb R}^{m\times n}$ for an arbitrary matrix, and
def\/ine operators   ${\cal A}_k: L^2(\Omega,{\mathbb R}^{n})\rightarrow L^2(\Omega,{\mathbb R}^{m})$ and
${\mathcal E}_{v_k v_k}, {\mathcal E}^{\dag}_{v_k v_k}, ({\mathcal E}^{1/2}_{v_k v_k})^{\dag} :L^{2}(\Omega,{\mathbb R}^{n})\rightarrow
L^{2}(\Omega,{\mathbb R}^{n})$ similarly to those in (\ref{qk}) and (\ref{ge}).

For the case under consideration (matrix ${\mathbb E}_{v_k v_k}$ is not invertible),
we introduce the SVD of  ${\mathbb E}_{xv_k}({\mathbb E}^{1/2}_{v_k v_k})^{\dag}$,
\begin{gather}
\label{svd2}
 U_k \Sigma_k V^{T}_k = {\mathbb E}_{xv_k}({\mathbb E}^{1/2}_{v_k v_k})^{\dag},
\end{gather}
 where, as above,
$U_k\in {\mathbb R}^{m\times n}$, $V_k\in {\mathbb R}^{n\times n}$
are  orthogonal and  $\Sigma_k\in {\mathbb R}^{n\times n}$ is diagonal,
\begin{gather}
\label{sqd3}
U_k = [s_{k1},\ldots,s_{kn}],\qquad  V_k =[d_{k1},\ldots,d_{kn}] \qquad\mbox{and}
\qquad \Sigma_k = \mbox{diag}(\beta_{k1},\ldots,\beta_{kn})
\end{gather}
 with $\beta_{k1} \geq
\cdots \geq \beta_{kr} > 0$, $\beta_{k,r+1} = \cdots = \beta_{kn} = 0$, $r=1,\ldots,n$ and
$r=r(k)$.

Let us set
\begin{gather}
U_{k \eta_k} = [s_{k1},\ldots,s_{k\eta_k}],\qquad V_{k\eta_k} = [d_{k1},\ldots, d_{k\eta_k}]\qquad
\mbox{and}\nonumber\\
\Sigma_{k\eta_k} = \mbox{diag}\,(\beta_{k1},\ldots,\beta_{k\eta_k}),\label{tr-svd}
\end{gather}
where $U_{k\eta_k}\in {\mathbb R}^{m\times \eta_k}$, $V_{k\eta_k}\in
{\mathbb R}^{n\times \eta_k}$  and  $\Sigma_{k\eta_k}\in {\mathbb R}^{\eta_k \times \eta_k}$.
 Now we  def\/ine $G_{k\eta_k} \in{\mathbb R}^{m\times n}$ and $\g_{k\eta_k}:L^{2}
 (\Omega,{\mathbb R}^{n})\rightarrow L^{2}(\Omega,{\mathbb R}^{m})$ by
\begin{gather}
\label{trsvd2}
G_{k\eta_k}  = U_{k\eta_k}\Sigma_{k\eta_k}V_{k\eta_k}^{T}\qquad\mbox{and}\qquad [\g_{k\eta_k}(\w_k)](\omega) = G_{k\eta_k}[\w_k(\omega)],
\end{gather}
respectively, for any $\w_k\in L^{2}(\Omega,{\mathbb R}^{n})$.

As noted before, we write ${\cal I}$ for the identity operator.

\begin{theorem}
\label{sol2}
Let $\vv_1,\ldots,\vv_p$ be determined by Lemma {\rm \ref{ort3}}.
Then  $f^0$ and ${\mathcal F}_1^0, \ldots, {\mathcal F}_p^0$,
 satisfying \eqref{min1}--\eqref{con1}, are determined by
\begin{gather}
\label{sol-f02}
f^0 = E[\x] - \sum_{k=1}^p F^0_k E[\vv_k]
\end{gather}
and
\begin{gather}
\label{sol-f12}
{\mathcal F}^0_1  = \g_{1\eta_1}({\mathcal E}^{1/2}_{v_1 v_1})^{ \dag}
+ {\cal A}_1[{\cal I} - {\mathcal E}_{v_1 v_1}^{1/2}({\mathcal E}^{1/2}_{v_1 v_1})^{\dag}],\\
\cdots \cdots \cdots\cdots \cdots \cdots \cdots \cdots \cdots \cdots \cdots \cdots \cdots \cdots\nonumber \\
\label{sol-fp2}
{\mathcal F}^0_p = \g_{p\eta_p}({\mathcal E}^{1/2}_{v_p v_p})^{ \dag}
+ {\cal A}_p[{\cal I} - {\mathcal E}_{v_p v_p}^{1/2}({\mathcal E}^{1/2}_{v_p v_p})^{\dag}],
\end{gather}
where for $k=1, \ldots,p,$ ${\cal A}_k$ is any linear operator such that $\mbox{\rm rank}\,
 {\mathcal F}^0_k \leq \eta_k$\footnote{In particular,  ${\cal A}_k$ can be chosen as the zero operator.}.

The accuracy associated with transform  ${\mathcal T}_p^0$ given by  \eqref{th1}
and  \eqref{sol-f02}--\eqref{sol-fp2} is such that
\begin{gather}
\label{er12}
E[\|\x - {\mathcal T}_p^0(\y)\|^2] =\|{\mathbb E}_{xx}^{1/2}\|^2 -\sum_{k=1}^p \sum_{j=1}^{\eta_k} \beta^2_{kj}.
\end{gather}
\end{theorem}

\begin{proof}  For  $\vv_1,\ldots,\vv_p$ determined by Lemma \ref{ort3}, $J(f,{\mathcal F}_1, \ldots ,{\mathcal F}_p)$
is represented by (\ref{jqa}) as well.
Let us consider $J_0$,  $J_1$ and $J_2$ given by
\begin{gather}
\label{j02}
J_0 = \|{\mathbb E}_{xx}^{1/2}\|^2 -\sum_{k=1}^p  \|{\mathbb E}_{xv_k}({\mathbb E}_{v_k v_k}^{1/2})^{\dag} \|^2,
\\
\label{j122}
J_1 = \|f - E[\x] + \sum_{k=1}^p F_k E[\vv_k]\|^2\qquad\mbox{and}\qquad
J_2 =\sum_{k=1}^p \|F_k{\mathbb E}_{v_kv_k}^{1/2} - {\mathbb E}_{xv_k}({\mathbb E}_{v_kv_k}^{1/2})^{\dag}\|^2.
\end{gather}
To show that
\begin{gather}
\label{jff}
J(f,{\mathcal F}_1, \ldots ,{\mathcal F}_p) = J_0 + J_1 + J_2
\end{gather}
with $J(f,{\mathcal F}_1, \ldots ,{\mathcal F}_p)$ def\/ined by (\ref{jqa}),  we use the relationships (see \cite{tor1})
\begin{gather}
\label{eee1}
{\mathbb E}_{xv_k}{\mathbb E}_{v_k v_k}^\dag {\mathbb E}_{v_k v_k}
= {\mathbb E}_{xv_k} \qquad\mbox{and}\qquad {\mathbb E}_{v_k v_k}^{\dag}{\mathbb E}_{v_k v_k}^{ 1/2}
= ({\mathbb E}_{v_k v_k}^{1/2})^{\dag}
\end{gather}
Then
\begin{gather}
J_1 = \mbox{tr} \Bigg(ff^T - fE[\x^T] + \sum_{k=1}^p f E[\vv_k^T]F_k + E[\x] E[\x^T]
 -  E[\x]f^T - \sum_{k=1}^p E[\x] E[\vv_k^T] F_k^T\nonumber\\
 \phantom{J_1 =}{} +  \sum_{k=1}^p F_k E[\vv_k]f^T
- \sum_{k=1}^p F_k E[\vv_k] E[\x^T] + \sum_{k=1}^p F_k E[\vv_k]  \sum_{i=1}^p E[\vv_i^T] F_i^T\Bigg)\label{j112}
\end{gather}
and
\begin{gather}
J_2  =  \sum_{k=1}^p \mbox{tr}
(F_k - {\mathbb E}_{xv_k}{\mathbb E}_{v_kv_k}^\dag){\mathbb E}_{v_kv_k}(F_k^T
- {\mathbb E}_{v_kv_k}^\dag{\mathbb E}_{v_k x})\nonumber\\
\phantom{J_2 }{} = \sum_{k=1}^p \mbox{tr} (F_k {\mathbb E}_{v_kv_k} F_k^T
- F_k {\mathbb E}_{v_k x} - {\mathbb E}_{xv_k}F_k^T + {\mathbb E}_{xv_k}{\mathbb E}_{v_kv_k}^\dag {\mathbb E}_{v_kx}),
\label{j22}
\end{gather}
where
\begin{gather}
\label{aea2}
\sum_{k=1}^p \mbox{tr} (F_k {\mathbb E}_{v_kv_k} F_k^T)= \mbox{tr}\Bigg[E\Bigg(\sum_{k=1}^p F_k
\vv_k \sum_{i=1}^p \vv_i^T F_i^T\Bigg)\Bigg] -\mbox{tr}\Bigg(\sum_{k=1}^p F_k E[\vv_k]  \sum_{i=1}^p E[\vv_i^T] F_i^T\Bigg)
\end{gather}
because
\begin{gather}
\label{v-ort2}
E[\vv_i \vv_k^T] - E[\vv_i] E[ \vv_k^T] = {\mathbb O}\qquad\mbox{for}\quad  i\neq k
\end{gather}
 due to orthogonality of the  vectors $\vv_1, \ldots, \vv_s$.
On the basis of (\ref{eee1})--(\ref{aea2}) and  similarly to (\ref{jj012})--(\ref{nj012}),
we establish that (\ref{jff}) is true.
 Hence,
\begin{gather}
 J(f,{\mathcal F}_1, \ldots ,{\mathcal F}_p)  =
\|{\mathbb E}_{xx}^{1/2}\|^2 -
\sum_{k=1}^p \|{\mathbb E}_{xv_k}({\mathbb E}_{v_k v_k}^{1/2})^{\dag} \|^2
+ \|f - E[\x]  \nonumber  \\
\label{proof22}
\phantom{J(f,{\mathcal F}_1, \ldots ,{\mathcal F}_p)  =}{}
+ \sum_{k=1}^p F_k E[\vv_k]\|^2+   \sum_{k=1}^p \|F_k{\mathbb E}_{v_k v_k}^{1/2} - {\mathbb E}_{xv_k}({\mathbb E}_{v_k v_k}^{1/2})^{\dag}\|^2.
\end{gather}
It follows from the last two terms in (\ref{proof22}) that
the constrained minimum (\ref{min1})--(\ref{con1}) is achieved
if $f=f^0$ with $f^0$ given by  (\ref{sol-f02}),  and  $F_k^0$ is such that
\begin{gather}
\label{min-k2}
J_k(F_k^0) = \min_{F_k} J_k(F_k) \qquad \mbox{subject to} \quad \mbox{rank} \,(F_k) = \eta_k,
\end{gather}
where $J_k(F_k) = \|F_k{\mathbb E}_{v_k v_k}^{1/2} - {\mathbb E}_{xv_k}({\mathbb E}_{v_k v_k}^{1/2})^{\dag}\|^2$.
The constrained minimum (\ref{min1})--(\ref{con1})
is achieved if $f=f^0$ is def\/ined by  (\ref{sol-f02}), and if \cite{gol1}
\begin{gather}
\label{me}
F_k{\mathbb E}_{v_k v_k}^{1/2} = G_{\eta_k}.
\end{gather}
The matrix equation (\ref{me}) has the general solution  \cite{ben1}
\begin{gather}
\label{proof312}
F_k = F^0_k =G_{\eta_k}({\mathbb E}_{v_k v_k}^{1/2})^\dag
+ A_k[I - {\mathbb E}_{v_k v_k}^{1/2} ({\mathbb E}_{v_k v_k}^{1/2})^\dag]
\end{gather}
if and only if
\begin{gather}
\label{proof42}
G_{\eta_k}({\mathbb E}_{v_k v_k}^{1/2})^\dag {\mathbb E}_{v_k v_k}^{1/2}= G_{\eta_k}.
\end{gather}
 The latter is satisf\/ied on the basis of the following
 derivation\footnote{Note that the matrix
$I - {\mathbb E}_{v_k v_k}^{1/2} ({\mathbb E}_{v_k v_k}^{1/2})^\dag$
is simply a projection onto the null space of ${\mathbb E}_{v_k v_k}$
 and can be replaced by $I - {\mathbb E}_{v_k v_k} ({\mathbb E}_{v_k v_k})^\dag$.}.

As an extension of the technique presented in the proving Lemmata~1
and~2 in \cite{tor1}, it can be shown that for any matrices $Q_1, Q_2\in{\mathbb R}^{m\times n}$,
\begin{gather}
\label{proof52}
{\cal N}(Q_1)\subseteq  {\cal N}(Q_2)\quad \Rightarrow \quad Q_2(I- Q_1^\dag Q_1) = \mathbb O,
\end{gather}
 where ${\cal N}(Q_i)$ is the null space of $Q_i$ for $i=1,2$. In regard of the equation under consideration,
\begin{gather}
\label{proof62}
{\cal N}([{\mathbb E}_{v_k v_k}^{1/2}]^\dag)
\subseteq  {\cal N}({\mathbb E}_{x v_k}[{\mathbb E}_{v_k v_k}^{1/2}]^\dag).
\end{gather}
 The def\/inition of $G_{\eta_k}$ implies that
\[
{\cal N}({\mathbb E}_{x v_k}[{\mathbb E}_{v_k v_k}^{1/2}]^\dag)
\subseteq {\cal N}(G_{\eta_k})\qquad\mbox{and then}\qquad {\cal N}([{\mathbb E}_{v_k v_k}^{1/2}]^\dag)
 \subseteq {\cal N}(G_{\eta_k}).
\]
 On the basis of (\ref{proof52}), the latter implies
 $G_{\eta_k}[I-({\mathbb E}_{v_k v_k}^{1/2})^\dag {\mathbb E}_{v_k v_k}^{1/2}]={\mathbb O}$,
 i.e.\ (\ref{proof42}) is true. Hence,  (\ref{proof312}) and (\ref{sol-f12})--(\ref{sol-fp2}) are true as well.

Next,  similar to (\ref{kr}),
\begin{gather}
\label{g-e}
 \|{\mathbb E}_{xv_k}({\mathbb E}_{v_kv_k}^{1/2})^\dag\|^2 -\|G_{\eta_k}
  - {\mathbb E}_{xv_k}({\mathbb E}_{v_kv_k}^{1/2})^\dag\|^2 =
  \sum_{j=1}^{\eta_k} \beta^2_{kj}.
\end{gather}
Then (\ref{er12}) follows from (\ref{proof22}), (\ref{proof312}), (\ref{sol-f02}) and (\ref{g-e}).
\end{proof}

\begin{remark}\label{remark4}
The known reduced-rank transforms   based on the Volterra polynomial structure
\cite{yam2,tor3,tor4} require the computation of a covariance matrix similar
to  ${\mathbb E}_{vv}$, where $\vv = [\vv_1,\ldots,\vv_p]^T$,
but for $p=N$ where $N$ is large (see Sections \ref{intr} and \ref{summ}).
The relationships (\ref{jj012})--(\ref{min-k}) and (\ref{j112})--(\ref{min-k2})
illustrate the nature of the proposed method and its dif\/ference from the
techniques in \cite{yam2,tor3,tor4}: due to the structure  (\ref{t3})
of the transform ${\mathcal T}_p$,  the procedure for f\/inding $f^0$,
${\mathcal F}_1^0$, $\ldots$,  ${\mathcal F}_p^0$ avoids  direct computation of ${\mathbb E}_{vv}$
which could be  troublesome due to large $N$. If operators  ${\mathcal Q}_1,\ldots,{\mathcal Q}_p$
are orthonormal, as in Theorem~\ref{sol1}, then (\ref{v-ort}) is true and
the covariance matrix ${\mathbb E}_{vv}$ is reduced to the identity.
If operators  ${\mathcal Q}_1,\ldots,{\mathcal Q}_p$ are orthogonal,
as in Theorem \ref{sol2}, then (\ref{v-ort2}) holds and  the covariance matrix ${\mathbb E}_{vv}$
is reduced to a block-diagonal form with non-zero blocks  ${\mathbb E}_{v_1v_1}, \ldots, {\mathbb E}_{v_pv_p}$
so that
\[
 {\mathbb E}_{vv}=\left [ \begin{array}{cccc}
{\mathbb E}_{v_1v_1} &    {\mathbb O}   & \ldots & {\mathbb O} \\
{\mathbb O}      & {\mathbb E}_{v_2v_2} & \ldots &  {\mathbb O} \\
\ldots & \ldots & \ldots & \ldots \\
{\mathbb O}      & {\mathbb O}   & \ldots &  {\mathbb E}_{v_p v_p}
\end{array} \right ]
\]
 with ${\mathbb O}$ denoting the zero block. As a result,
 the procedure for f\/inding  $f^0$, ${\mathcal F}_1^0, \ldots,  {\mathcal F}_p^0$
 is reduced to  $p$ separate rank-constrained problems (\ref{min-k}) or (\ref{min-k2}).
 Unlike the methods in \cite{yam2,tor3,tor4},  the operators ${\mathcal F}^0_1,\ldots,{\mathcal F}^p_0$
 are determined with {\em  much smaller} $m\times n$ and {$n\times n$ matrices}
 given by the simple formulae (\ref{sol-f0}) and (\ref{sol-f02})--(\ref{sol-fp2}).
 This implies a reduction in computational work compared with that required by the approach
 in \cite{tor3,tor4,how2}.
\end{remark}

\begin{corollary}\label{corollary2} Let $\vv_1,\ldots,\vv_p$ be determined by Lemma {\rm\ref{ort3}}.
Then the vector $\bar{f}$ and operators $\bar{\mathcal F}_1, \ldots, \bar{\mathcal F}_p$,
satisfying the unconstrained minimum \eqref{min1}, are determined by
\begin{gather}
\label{sol-fc0}
\bar{f} = E[\x] - \sum_{k=1}^p \bar{F}_k E[\vv_k]
\end{gather}
and
\begin{gather}
\label{sol-fc2}
\bar{\mathcal F}_1  = {\mathcal E}_{xv_1}{\mathcal E}_{v_1 v_1}^{ \dag}
+ {\cal A}_1[{\cal I} - {\mathcal E}_{v_1 v_1}{\mathcal E}_{v_1 v_1}^{\dag}],\\
\cdots \cdots\cdots \cdots\cdots \cdots\cdots \cdots\cdots \cdots\cdots \cdots \nonumber \\
\label{sol-fpc2}
\bar{\mathcal F}_p  = {\mathcal E}_{xv_p}{\mathcal E}_{v_p v_p}^{ \dag}
+ {\cal A}_p[{\cal I} - {\mathcal E}_{v_p v_p}{\mathcal E}_{v_p v_p}^{\dag}].
\end{gather}
The associated accuracy for transform  $\bar{{\mathcal T}}_p$, defined by
\[
\bar{{\mathcal T}}_p(\y) = \bar{f} + \sum _{k=1}^p \bar{\mathcal F}_k(\vv_k),
\]
 is given by
\begin{gather}
\label{er-c2}
E[\|\x - \bar{{\mathcal T}}_p(\y)\|^2] =\|{\mathbb E}_{xx}^{1/2}\|^2
-\sum_{k=1}^p  \|{\mathbb E}_{xv_k}({\mathbb E}^{1/2}_{v_k v_k})^{ \dag}\|^2.
\end{gather}
\end{corollary}

\begin{proof} It follows from (\ref{proof22}) that the unconstrained minimum (\ref{min1})
is achieved if $f$ is def\/ined by~(\ref{sol-fc0}) and if ${F}_k$ satisf\/ies the equation
$
F_k{\mathbb E}_{v_k v_k}^{1/2} - {\mathbb E}_{xv_k}({\mathbb E}_{v_k v_k}^{1/2})^{\dag} = {\mathbb O}
$
for each $k=1,\ldots,p$.
Similar to (\ref{me})--(\ref{proof312}), its general solution is given by
\[
F_k = \bar{F}_k ={\mathbb E}_{xv_k}{\mathbb E}_{v_k v_k}^{\dag}
+ A_k[I - {\mathbb E}_{v_k v_k}{\mathbb E}_{v_k v_k}^{\dag}]
\]
because ${\mathbb E}_{v_k v_k}^{1/2}({\mathbb E}_{v_k v_k}^{1/2})^{\dag}
 ={\mathbb E}_{v_k v_k} {\mathbb E}_{v_k v_k}^{\dag}$.
We def\/ine $\bar{\mathcal F}_k$  by $[\bar{\mathcal F}_k(\w_k)](\omega)
= \bar{F}_k[\w_k(\omega)]$  for all $k=1,\ldots,p$, and then (\ref{sol-fc2})--(\ref{sol-fpc2})
are true. The relation (\ref{er-c2}) follows from (\ref{proof22}) and  (\ref{sol-fc0})--(\ref{sol-fpc2}).
\end{proof}

\begin{remark} The dif\/ference between the transforms given by Theorems \ref{sol1} and \ref{sol2} is that
${\mathcal F}^0_k$ by~(\ref{sol-f0}) (Theorem~\ref{sol1}) does not contain a factor
associated with $({\mathbb E}_{v_kv_k}^{1/2})^\dag$ for all $k=1,\ldots.p$.
A similar observation is true for Corollaries \ref{corollary1} and \ref{corollary2}.
\end{remark}

\begin{remark}
The transforms given by Theorems \ref{sol1} and \ref{sol2} are not unique due to arbitrary operators
${\mathcal A}_1, \ldots, {\mathcal A}_p$. A natural particular choice
is  ${\mathcal A}_1 = \cdots = {\mathcal A}_p =  \mathbb O$.
\end{remark}

\subsubsection[Compression procedure by ${\mathcal T}^0_p$]{Compression procedure by $\boldsymbol{{\mathcal T}^0_p}$}
\label{sec5.2.3}

Let us consider transform ${\mathcal T}^0_p$ given by (\ref{th1}),
(\ref{sol-f02})--(\ref{sol-fp2}) with $A_k={\mathbb O}$ for $k=1,\ldots,p$ where $A_k$
is the matrix given in (\ref{proof312}).
We write $[{\mathcal T}^0_p(\y)](\omega) = T^0_p(y)$
with $T^0_p:{\mathbb R}^n\rightarrow {\mathbb R}^m$.

Let
\[
B^{(1)}_k = S_{k\eta_k}V_{k\eta_k}D_{k\eta_k}^{T}\qquad\mbox{and}
\qquad B^{(2)}_k = D_{k\eta_k}^{T}({\mathbb E}^{1/2}_{v_k v_k})^{\dag}
\]
so that $B^{(1)}_k\in{\mathbb R}^{m\times \eta_k}$
and $B^{(2)}_k\in{\mathbb R}^{\eta_k\times n}$. Here, $\eta_1$, $\ldots$, $\eta_p$
are determined by (\ref{con1}). Then
\[
T^0_p(y) = f + \sum_{k=1}^p B^{(1)}_k B^{(2)}_k v_k,
\]
where $v_k = \vv_k(\omega)$ and $B^{(2)}_k v_k\in{\mathbb R}^{\eta_k}$ for $k=1,\ldots,p$
with $\eta_1 + \cdots + \eta_p <m$.
Hence, matrices $B^{(2)}_1,\ldots,B^{(2)}_p$ perform compression
of the data presented by $v_1, \ldots, v_p$. Matrices $B^{(1)}_1,\ldots,B^{(1)}_p$
perform  reconstruction of the reference signal from the compressed data.

 The compression ratio of transform ${\mathcal T}^0_p$ is given by
\begin{gather}
\label{cr}
r^0 = (\eta_1 + \cdots + \eta_p)/m.
\end{gather}

\subsubsection[A special  case of transform ${\mathcal T}_p$]{A special  case of transform $\boldsymbol{{\mathcal T}_p}$}
\label{spec}

The results above have been derived for  any operators ${\boldsymbol{\varphi}}_{1},
\ldots,{\boldsymbol{\varphi}}_{p}$ in the model ${\mathcal T}_p$.
Some  specializations for ${\boldsymbol{\varphi}}_{1},\ldots,{\boldsymbol{\varphi}}_{p}$
were given in Section~\ref{some}. Here and in Section~\ref{part-cas},
we consider alternative forms for ${\boldsymbol{\varphi}}_{1},\ldots,{\boldsymbol{\varphi}}_{p}$.

 {\bf (i)}  Operators ${\boldsymbol{\varphi}}_{1},\ldots,{\boldsymbol{\varphi}}_{p}$ can be determined by
 a recursive procedure given below. The motivation follows
 from the observation that  performance of the transform ${\mathcal T}_p$
 is improved if $\y$ in (\ref{j1}) is replaced by an estimate of $\x$.

First, we set ${\boldsymbol{\varphi}}_{k}(\y) = \y$ and determine
estimate $\x^{(1)}$ of $\x$ from the solution of
problem~(\ref{min1}) (with no constraints (\ref{con1})) by
Corollaries \ref{corollary1} or \ref{corollary2} with $p = 1$.
Next, we put  ${\boldsymbol{\varphi}}_{1}(\y) = \y$ and
${\boldsymbol{\varphi}}_{2}(\y) =\x^{(1)}$, and f\/ind  estimate
$\x^{(2)}$ from the solution of unconstrained problem
(\ref{min1})  with $p = 2$. In general, for $j=1,\ldots,p$, we
def\/ine
 ${\boldsymbol{\varphi}}_{j}(\y)$ $=\x^{(j-1)}$, where $\x^{(j-1)}$ has been
 determined  similarly to  $\x^{(2)}$ from the previous steps. In particular,  $\x^{(0)} = \y$.

{\bf (ii)} Operators ${\boldsymbol{\varphi}}_1,\ldots,{\boldsymbol{\varphi}}_p$
can also be chosen as elementary functions. An example is given in item {\bf (i)}
of Section \ref{some} where  ${\boldsymbol{\varphi}}_k(\y)$ was constructed from
the power functions. An alternative possibility is to choose trigonometric
functions for constructing  ${\boldsymbol{\varphi}}_k(\y)$. For instance, one can put
\begin{gather}
\label{cos}
[{\boldsymbol{\varphi}}_1(\y)](\omega) =y \qquad \mbox{and}\qquad [{\boldsymbol{\varphi}}_{k+1}(\y)](\omega)
=  [\cos(k y_1), \ldots, \cos(k y_n)]^T
\end{gather}
with  $y= [y_1,\ldots,y_n]^T$ and $k=1,\ldots,p-1$.  In this paper, we do not analyse
such a possible choice for ${\boldsymbol{\varphi}}_1,\ldots,{\boldsymbol{\varphi}}_p$.

\subsubsection[Other particular cases of transform ${\mathcal T}_p$ and comparison with
known transforms]{Other particular cases of transform $\boldsymbol{{\mathcal T}_p}$\\ and comparison with
known transforms}
\label{part-cas}

 {\bf (i)} {\bf\em Optimal non-linear filtering}.
The transforms  $\hat{{\mathcal T}}_p$ (\ref{sol-q})--(\ref{fr1}) and  $\bar{{\mathcal T}}_p$
(\ref{sol-fc0})--(\ref{sol-fpc2}), which are  particular cases of the transforms given in
Theorems \ref{sol1} and \ref{sol2},
represent optimal f\/ilters   that perform pure f\/iltering with no signal compression. Therefore they
are important in their own right.

{\bf (ii)} {\bf\em The Fourier series as a particular case of  transform $\bar{{\mathcal T}}_p$.}
 For the case of the minimization problem (\ref{min1}) with no constraint (\ref{con1}),
 ${\mathcal F}_1, \ldots, {\mathcal F}_p$ are determined by the expressions
 (\ref{sol-q}) and (\ref{sol-fc0})--(\ref{sol-fpc2}) which are similar
 to those for the Fourier coef\/f\/icients  \cite{cot1}.
 The structure of the model $ {\mathcal T}_{p}$ presented by (\ref{t3})
 is dif\/ferent, of course,  from that for the  Fourier series and Fourier polynomial
 (i.e.\  a truncated Fourier series) in Hilbert space \cite{cot1}.
 The dif\/ferences are that ${\mathcal T}_{p}$ transforms $\y$ (not $\x$ as
 the Fourier polynomial does) and that $ {\mathcal T}_{p}$ consists of a combination of three
  operators  ${\boldsymbol{\varphi}}_k$, ${\mathcal Q}_k$
and  ${\mathcal F}_k$ where ${\mathcal F}_k:L^2(\Omega, \tilde{H}_k)\rightarrow L^2(\Omega, H_X)$
is an operator, not a scalar as in the Fourier series~\cite{cot1}.
The solutions (\ref{sol-q}) and (\ref{sol-fc0})--(\ref{sol-fpc2})
of the unconstrained problem (\ref{min1}) are given in terms of
the observed vector $\y$, not in terms of the basis of $\x$  as in the Fourier series/polynomial.
The special features of ${\mathcal T}_{p}$ require special computation methods as described in Section~\ref{sol}.

Here, we show that the  Fourier series is a particular case of the transform ${\mathcal T}_p$.

 Let $\x\in L^2(\Omega, H)$ with $H$ a Hilbert space, and
 let $\{ \vv_1, \vv_2, \ldots\}$ be an orthonormal basis in $L^2(\Omega, H).$
For any $\gf, \h \in  L^2(\Omega, H),$
 we def\/ine the scalar product $\langle \cdot, \cdot \rangle$ and the norm $\|\cdot\|_{_E}$ in $L^2(\Omega, H)$  by
\begin{gather}
\label{scal}
\langle \gf, \h \rangle = \int_{\Omega} \gf(\omega) \h(\omega)d\mu (\omega)\qquad\mbox{and}\qquad
\|\gf\|_{_E} = \langle \gf, \gf \rangle^{1/2},
\end{gather}
respectively.
{\samepage In particular, if  $H = {\mathbb R}^m$ then
\begin{gather}
\label{hilb}
\|\gf\|^2_{_E} = \int_{\Omega} \gf(\omega) [\gf(\omega)]^Td\mu (\omega)
 = \int_{\Omega} \|\gf(\omega)\|^2 d\mu (\omega)=  E[\|\gf\|^2],
\end{gather}
i.e.\ $E[\|\gf\|^2]$ is def\/ined similarly to that in  (\ref{e11}).}

Let us  consider the special case of transform ${\mathcal T}_p$
presented in item {\bf (iii)} of Section \ref{some}  and let
us also consider the unconstrained problem (\ref{min1}) formulated
in terms of such a ${\mathcal T}_p$ where we now assume that $\x$
has the zero mean, $f={\mathbb O}$, $p = \infty$, $\{\vv_1, \vv_2,\ldots \}$
is an orthonormal basis in $L^2(\Omega, H)$ and ${\mathcal F}_k$
is a scalar, not an operator as before. We denote $ \alpha_k = {\mathcal F}_k$
with $\alpha_k \in {\mathbb R}$. Then  similar to  (\ref{sol-q})
in Corollary \ref{corollary1}, the solution to unconstrained problem (\ref{min1}) is def\/ined by  $\hat{\alpha}_k$ such that
\[
\hat{\alpha}_k =  {\mathbb E}_{x v_k}  \qquad\mbox{with}\quad k=1,2,\ldots .
\]
Here, ${\mathbb E}_{x v_k} = E[\x \vv_k] - E[\x] E[ \vv_k] = E[\x \vv_k] = \langle \x, \vv_k \rangle $
 since $E[\x] = 0$ by the assumption.
Hence, $\hat{\alpha}_k = {\mathbb E}_{x v_k}$
is the Fourier coef\/f\/icient and the considered particular case of ${\mathcal T}_p(\y)$
with ${\mathcal F}_k$ determined by $\hat{\alpha}_k$ is given by
\begin{gather}
\label{f-ser}
{\mathcal T}_p(\y) = \sum_{k=1}^\infty \langle \x, \vv_k \rangle \vv_k.
\end{gather}
Thus, the Fourier series (\ref{f-ser}) in Hilbert space follows from (\ref{t3}), (\ref{min1})
 and (\ref{sol-q}) when ${\mathcal T}_p$ has the form given in  item {\bf (iii)} of Section \ref{some}
 with $\x$, $f$, $p$,  $\{\vv_1, \vv_2,\ldots \}$  and ${\mathcal F}_k$ as above.

{\bf (iii)} {\bf\em The  Wiener filter as a particular case of transform} $\bar{{\mathcal T}}_p$
(\ref{sol-fc0})--(\ref{sol-fpc2}). In the following Corollaries \ref{corollary3} and \ref{corollary4}, we show that
the f\/ilter  $\bar{{\mathcal T}}_p$  guarantees  better accuracy than that of the Wiener f\/ilter.

\begin{corollary}\label{corollary3} Let $p=1$, $E[\x]=0$,  $E[\y]=0$,
${\boldsymbol{\varphi}}_1={\cal I}$, ${\mathcal Q}_1={\cal I}$ and
$A_1={\mathbb O}$ or $A_1=E_{xy}E_{yy}^{\dag}$. Then $\bar{{\mathcal T}}_p$ is reduced to the filter
$\check{{\mathcal T}}$ such that
\[
[\check{{\mathcal T}}(\y)](\omega) = \check{T}[\y(\omega)]
\]
 with
\begin{gather}
\label{win2}
\check{T} = E_{xy}E_{yy}^\dag.
\end{gather}
\end{corollary}

\begin{remark}
The  unconstrained linear f\/ilter, given by (\ref{win2}), has been proposed in \cite{hua1}.
The f\/ilter~(\ref{win2}) is treated as a generalisation of the Wiener f\/ilter.
\end{remark}

Let  $\tilde{\x}$, $\tilde{\vv}_1,\ldots,\tilde{\vv}_p$ be the
zero mean vectors. The transform $\bar{{\mathcal T}}_p$,
applied to $\tilde{\x}$, $\tilde{\vv}_1,\ldots,\tilde{\vv}_p$, is denoted by $\bar{{\mathcal T}}_{W, p}$.

\begin{corollary}\label{corollary4} The error $E[\|\tilde{\x} - \bar{{\mathcal T}}_{W, p}(\tilde{\y})\|^2]$
associated with the transform $\bar{{\mathcal T}}_{W, p}$
is smaller than the error $E[\|\tilde{\x} - \check{{\mathcal T}}(\tilde{\y})\|^2] $
 associated with the  Wiener filter {\rm \cite{hua1}} by $\sum\limits_{k=2}^p
  \|E_{\tilde{x}\tilde{v}_k}(E^{1/2}_{\tilde{v}_k \tilde{v}_k})^{\dag}\|^2$, i.e.
\begin{gather}
\label{er-wien1}
E[\|\tilde{\x} - \bar{{\mathcal T}}_{W, p}(\tilde{\y})\|^2]
=E[\|\tilde{\x} - \check{{\mathcal T}}(\tilde{\y})\|^2]
-\sum_{k=2}^p  \|E_{\tilde{x}\tilde{v}_k}(E^{1/2}_{\tilde{v}_k \tilde{v}_k})^{ \dag}\|^2.
\end{gather}
\end{corollary}

\begin{proof} It is easy to show that
\begin{gather}
\label{er-wien2}
E[\|\tilde{\x} - \check{{\mathcal T}}(\tilde{\y})\|^2]
= \|E_{\tilde{x}\tilde{x}}^{1/2}\|^2 - \|E_{\tilde{x}\tilde{v_1}}(E^{1/2}_{\tilde{v_1} \tilde{v_1}})^{ \dag}\|^2,
\end{gather}
and then (\ref{er-wien1}) follows from (\ref{er-c2}) and (\ref{er-wien2}).
\end{proof}

{\bf (iv)} {\bf\em  The KLT as a particular case of transform} ${\mathcal T}_p^0$  (\ref{sol-f02})--(\ref{er12}).
 The KLT \cite{hua1} follows from  (\ref{sol-f02})--(\ref{er12}) as a particular case if
$f={\mathbb O}$, $p=1$, ${\boldsymbol{\varphi}}_1={\mathcal I}$,   ${\mathcal Q}_1={\mathcal I}$
and $A_1={\mathbb O}$.

To compare the transform ${\mathcal T}^0_p$ with the KLT \cite{hua1},
we apply ${\mathcal T}^0_p$, represented by  (\ref{sol-f02})--(\ref{er12}),  to
the zero mean vectors $\tilde{\x}$, $\tilde{\vv}_1,\ldots,\tilde{\vv}_p$
as above. We write ${\mathcal T}^*_p$ for such a version
of ${\mathcal T}^0_p$, and ${\mathcal T}_{_{\rm KLT}}$ for the KLT  \cite{hua1}.

\begin{corollary}\label{corollary5}
The error $E[\|\tilde{\x} - {{\mathcal T}}^*_p(\tilde{\y})\|^2]$
associated with the transform ${\mathcal T}^*_p$ is smaller
than the error $E[\|\tilde{\x} - {{\mathcal T}}_{_{\rm KLT}}(\tilde{\y})\|^2]$
associated with the KLT {\rm \cite{hua1}} by $\sum\limits_{k=2}^p \sum\limits_{j=1}^{\eta_k} \beta^2_{kj}$, i.e.
\begin{gather}
\label{er-c3}
E[\|\tilde{\x} - {{\mathcal T}}^*_p(\tilde{\y})\|^2]
= E[\|\tilde{\x} - {{\mathcal T}}_{_{\rm KLT}}(\tilde{\y})\|^2]  - \sum_{k=2}^p \sum_{j=1}^{\eta_k} \beta^2_{kj}.
\end{gather}
\end{corollary}

\begin{proof} The error associated with ${\mathcal F}_{_{\rm KLT}}$
\cite{hua1} is represented by (\ref{er12}) for $p=1$,
\begin{gather}
\label{er-klt}
E[\|\tilde{\x} - {{\mathcal T}}_{_{\rm KLT}}(\tilde{\y})\|^2]
= \|E_{\tilde{x}\tilde{x}}^{1/2}\|^2 - \sum_{j=1}^{\eta_1} \beta^2_{1j}.
\end{gather}
Then (\ref{er-c3}) follows from (\ref{er12}) and (\ref{er-klt}).
\end{proof}

{\bf (v)} {\bf\em The transform {\rm \cite{tor1}} as a particular case of transform ${\mathcal T}^0_p$.}
The transform \cite{tor1}  follows from (\ref{t3}) as a particular case
if  $f={\mathbb O}$, $p=2$, ${\boldsymbol{\varphi}}_1(\y)=\y$,
${\boldsymbol{\varphi}}_2(\y)=\y^2$ and ${\mathcal Q}_1={\mathcal Q}_2={\cal I}$
where $\y^2$ is def\/ined by $\y^2(\omega) = [y_1^2, \ldots, y_n^2]^T$.
We note that transform \cite{tor1} has been generalized in \cite{tor3}.

{\bf (vi)} {\bf\em The transforms {\rm \cite{tor3}}
 as  particular cases of transform ${\mathcal T}_p$.}
 The transform  \cite{tor3} follows from (\ref{t3}) if ${\mathcal Q}_k = {\cal I}$,
 ${\boldsymbol{\varphi}}_k(\y) = \y^k$ where
 $\y^k = (\y,\ldots, \y)\in L^{2}(\Omega,{\mathbb R}^{nk})$,
 ${\mathbb R}^{nk}$ is the $k$th degree of~${\mathbb R}^{n}$,
 and if ${\mathcal F}_k$ is a $k$-linear operator.

To compare  transform ${\mathcal T}^0_p$ and  transform  ${\mathcal T}_{\mbox{\scriptsize \cite{tor3}}}$
\cite{tor3} of rank $r$, we write $z_j = y_j y$, $z = [z_1,\ldots,z_n]^T$, $s = [1\hspace*{1mm}
y^T\hspace*{1mm}  z^T]^T$ and denote by $\alpha_1,\ldots,\alpha_r$ the non-zero singular values
associated with the truncated SVD for the matrix ${\mathbb E}_{xs} ({\mathbb E}_{ss}^{1/2})^\dag$. Such a SVD
is constructed similarly to that in (\ref{svd2})--(\ref{tr-svd}).

\begin{corollary}\label{corollary6}
Let $\Delta_p = \sum\limits_{k=1}^p \sum\limits_{j=1}^{\eta_k} \beta^2_{kj} - \sum\limits_{j=1}^{r} \alpha^2_{j}$
and let $\Delta_p \geq 0$.
The error $E[\|{\x} - {{\mathcal T}}^0_p({\y})\|^2]$ associated with the transform ${\mathcal T}^0_p$ is
less than the error $E[\|{\x} - {\mathcal T}_{\mbox{\rm \scriptsize \cite{tor3}}}(\y)\|^2]$ associated
with  the transform   ${\mathcal T}_{\mbox{\rm \scriptsize \cite{tor3}}}$ by $\Delta_p$, i.e.
\begin{gather}
\label{er-c20}
E[\|{\x} - {{\mathcal T}}^0_p({\y})\|^2]  = E[\|{\x} - {\mathcal T}_{\mbox{\rm \scriptsize \cite{tor3}}}(\y)\|^2]
- \Delta_p.
\end{gather}
\end{corollary}

\begin{proof} It follows from \cite{tor3} that
\begin{gather}
\label{er-20}
E[\|{\x} - {\mathcal T}_{\mbox{\rm \scriptsize \cite{tor3}}}({\y})\|^2]
= \|E_{{x}{x}}^{1/2}\|^2 - \sum_{j=1}^{r} \alpha^2_{j}.
\end{gather}
Then (\ref{er-c20}) follows from (\ref{er12}) and (\ref{er-20}).
\end{proof}

We note that, in general, a theoretical verif\/ication of the condition $\Delta_p \geq 0$ is not
straightforward. At the same time, for any  particular   $\x$ and $\y$, $\Delta_p $ can be
estimated numerically.

Although the transform ${{\mathcal T}}^0_p$ includes the transform
${\mathcal T}_{\mbox{\rm \scriptsize \cite{tor3}}}$, the accuracy of ${\mathcal T}_{\mbox{\rm \scriptsize \cite{tor3}}}$
is, in general, better than that of ${{\mathcal T}}^0_p$ for the same degrees of ${{\mathcal T}}^0_p$ and
${\mathcal T}_{\mbox{\rm \scriptsize \cite{tor3}}}$. This is because
${\mathcal T}_{\mbox{\rm \scriptsize \cite{tor3}}}$ implies more terms. For instance,
${\mathcal T}_{\mbox{\rm \scriptsize \cite{tor3}}}$ of degree two consists of $n+1$ terms while
${{\mathcal T}}^0_p$, for $p=2$, consists of three terms only. If for a given $p$, the condition  $\Delta_p \geq 0$
is not fulf\/illed, then the accuracy  $E[\|{\x} - {{\mathcal T}}^0_p({\y})\|^2]$ can be improved by
increasing $p$ or by applying the iterative method presented in \cite{tor2}.

{\bf (vii)} Unlike the techniques presented in \cite{ten1,row1},
our method implements simultaneous f\/iltering and compression,
and provides this data processing in probabilistic setting.
The idea of implicitly mapping the data into a high-dimensional feature
space \cite{vap1,scho1,cri1} could be extended to the transform presented in this paper.
We intend to develop such an extension in the  future.

\section{Numerical realization}
\label{num}

{\bf 6.1.} {\bf Orthogonalization.} Numerical
realization of transforms of random vectors implies a representation
of observed data and estimates of covariance matrices in the form of associated  samples.

For the random vector $\uu_k$, we have $q$ realizations, which are concatenated
into $ {n\times q}$ mat\-rix~$U_k$.  A  column of  $U_k$ is a realization
of  $\uu_k$. Thus, a sequence of vectors $\uu_1,\ldots, \uu_p$
is represented by a sequence of matrices $U_1,\ldots,U_p$.
Therefore the transformation of $\uu_1,\ldots, \uu_p$
to orthonormal or orthogonal vectors  $\vv_1,\ldots,\vv_p$
 (by Lemmata~\ref{lemma1} and~\ref{ort3}) is reduced to a procedure for
 matrices $U_1,\ldots,U_p$  and $V_1,\ldots,V_p$.
 Here, $V_k\in{\mathbb R}^{n\times q}$ is a matrix formed from
 realizations of the random vector $\vv_k$ for each $k=1,\ldots,p$.

Alternatively, matrices $V_1,\ldots,V_p$ can be determined
 from  known procedures for matrix orthogonalization \cite{gol1}.
In particular, the QR decomposition \cite{gol1}
can be exploited in the following way. Let us
form a matrix $U = [U_1^T \ldots U_p^T]^T\in{\mathbb R}^{np\times q}$
where $p$ and $q$ are chosen such that $np=q$, i.e.\ $U$
is square\footnote{Matrix $U$ can also be presented as $U = [U_1 \ldots U_p]$ with $p$ and $q$ such that $n=pq$.}.
 Let
\[
U = VR
\]
 be the QR decomposition for $U$ with $V\in{\mathbb R}^{np\times q}$
 orthogonal and $R\in{\mathbb R}^{np\times q}$ upper triangular.
\noindent
 Next, we write
$V =  [V_1^T \ldots V_p^T]^T\in{\mathbb R}^{mp\times q}$
where $V_k\in{\mathbb R}^{n\times q}$ for  $k=1,\ldots,p$.
The submatrices $V_1,\ldots,V_p$ of  $V$ are orthogonal, i.e.\
$
V_i V_j^T = \left \{ \begin{array}{@{}cc}
{\mathbb O}, &  i\neq j,\\
 I, &  i=j, \end{array}\right.$ for $i,j=1,\ldots,p,
$
as required.

Other known procedures for matrix orthogonalization can be applied to $U_1,\ldots,U_p$ in a similar fashion.

\begin{remark} For the cases when $\vv_1,\ldots,\vv_p$ are orthonormal
or orthogonal but not orthonormal, the associated  accuracies (\ref{er1}), (\ref{er2}), (\ref{er12}) and
(\ref{er-c2})   dif\/fer for the factors depending on $({\mathbb E}_{v_kv_k}^{1/2})^\dag$.
In the case of orthonormal $\vv_1,\ldots,\vv_p$, $({\mathbb E}_{v_kv_k}^{1/2})^\dag =I$
 and this circumstance can lead to an increase in the accuracy.
 \end{remark}

{\bf 6.2.} {\bf Covariance matrices.}
The expectations and covariance matrices in Lemmata \ref{lemma1}--\ref{ort3} and Theorems \ref{sol1}--\ref{sol2}
can  be estimated, for example,  by the techniques developed in  \cite{per1,kau1,sch3,kub1,led1,leu1}.
We note that such estimation procedures represent  specif\/ic problems which are not considered here.

{\bf 6.3.} {\bf ${{\mathcal T}^0_p}$, ${{\hat{{\mathcal T}}_p}}$
and ${{\bar{{\mathcal T}}_p}}$ for zero mean vectors.}
The computational work for ${\mathcal T}^0_p$ (Theorems~\ref{sol1} and~\ref{sol2}),
$\hat{{\mathcal T}}_p$ and $\bar{{\mathcal T}}_p$ (Corollaries~\ref{corollary1} and~\ref{corollary2})
can be reduced if ${\mathcal T}^0_p$, $\hat{{\mathcal T}}_p$ and $\bar{{\mathcal T}}_p$
are   applied to the zero mean vectors $\tilde{\x}$, $\tilde{\vv}_1, \ldots, \tilde{\vv}_p$
given by  $\tilde{\x} = \x - E[\x]$, $\tilde{\vv}_1 = \vv_1 - E[\vv_1],\ldots,\tilde{\vv}_p = \vv_p - E[\vv_p]$.
 Then $f^0={\mathbb O}$ and $\bar{f}={\mathbb O}$. The estimates of the original $\x$ are  then given by
\[
\check{\x} = E[\x] +\sum _{k=1}^p {\mathcal F}_k^0(\tilde{\vv}_k),
\qquad  \hat{\x} = E[\x] + \sum _{k=1}^p \hat{\mathcal F}_k(\tilde{\vv}_k)
\qquad\mbox{and}\qquad \bar{\x} = E[\x] + \sum _{k=1}^p \bar{\mathcal F}_k(\tilde{\vv}_k)
\]
respectively. Here, ${\mathcal F}_k^0$, $\hat{\mathcal F}_k$
and $\bar{\mathcal F}_k$ are def\/ined similarly to (\ref{sol-f0}),
(\ref{sol-q}), (\ref{sol-f12}), (\ref{sol-fp2}), (\ref{sol-fc2}) and (\ref{sol-fpc2}).

\section{Discussion}

Some distinctive features of the proposed techniques are summarized as follows.

\begin{remark}
It follows from Theorems \ref{sol1} and \ref{sol2},
and Corollaries \ref{corollary1} and \ref{corollary2} that the accuracy associated with the proposed transform improvs when  $p$  increases.
\end{remark}

\begin{remark} Unlike the approaches based on Volterra
polynomials \cite{tor3,tor4,how2} our method does not require
computation of pseudo-inverses for large $N\times N$ matrices
with $N=n+n^2+\cdots +n^{p-1}$. Instead, the proposed transforms
use pseudo-inverses of  $n\times n$  matrix ${\mathbb E}_{v_kv_k}$. See Theorems~\ref{sol1} and \ref{sol2}.
This  leads to a substantial reduction in computational work.
\end{remark}

\begin{remark}
The idea of the recurrent transform \cite{tor2} can
be extended for the proposed transform in a way similar to
that considered in  \cite{tor2}. 
The authors intend to develop a theory for such an extension in a feasible future.
\end{remark}

\section{Conclusions}

The new results obtained in the paper are summarized as follows.

We have proposed a new approach to constructing optimal nonlinear transforms for random vectors.
The approach is based on a representation of a transform in the form
of the sum of $p$ reduced-rank transforms. Each particular  transform
is formed by the linear reduced-rank operator ${\mathcal F}_k$, and by
operators ${\boldsymbol{\varphi}}_k$  and ${\mathcal Q}_k$ with $k=1,\ldots,p$.
Such a device  allows us to improve the numerical characteristics
(accuracy, compression ration and computational work) of the known transforms
 based on the Volterra polynomial structure \cite{tor3,tor4,how2}.
These objectives are achieved due to the special ``intermediate''
operators ${\boldsymbol{\varphi}}_1, \ldots,{\boldsymbol{\varphi}}_p$
 and ${\mathcal Q}_1, \ldots, {\mathcal Q}_p$.
In particular, we have proposed two types of orthogonalizing operators  ${\mathcal Q}_1, \ldots, {\mathcal Q}_p$
(Lemmata~\ref{lemma1} and~\ref{ort3}) and a specif\/ic method
for determining ${\boldsymbol{\varphi}}_1, \ldots, {\boldsymbol{\varphi}}_p$ (Section~\ref{spec}).
Such operators  reduce  the determination of optimal  linear reduced-rank
operators  ${\mathcal F}^0_1, \dots, {\mathcal F}^0_p$
to the computation of a sequence of relatively small matrices (Theorems~\ref{sol1} and \ref{sol2}).

Particular cases of the proposed transform, which follow from the
solution of the unconstrained minimization problem (\ref{min1}),
have been presented in Corollaries \ref{corollary1} and \ref{corollary2}. Such  transforms are
treated as  new optimal nonlinear f\/ilters and, therefore, are important in their own right.

The explicit representations of the accuracy associated with the
proposed transforms have been rigorously justif\/ied in Theorems \ref{sol1} and \ref{sol2}, and
Corollaries \ref{corollary1} and \ref{corollary2}.

It has been shown that the proposed approach generalizes the Fourier series
in Hilbert space (Section \ref{part-cas}), the Wiener f\/ilter,
the Karhunen--Lo\`{e}ve transform (KLT) and the  known optimal
transforms \cite{tor1,tor3,tor4}. See Corollaries \ref{corollary3}, \ref{corollary4}
and \ref{corollary5},  and Section \ref{part-cas} in this regard.
In particular, it has been shown  that the accuracies associated
with the proposed transforms are better  than those of the Wiener f\/ilter (Corollary~\ref{corollary4})
and the  KLT (Corollary \ref{corollary5}).

\appendix

\section{Appendix}

{\bf Proof of Lemma \ref{lemma1}.}  Let us write
\[
\w_1 =\uu_1\qquad\mbox{and}\qquad
\w_i = \uu_i - \sum_{k=1}^{i-1} {\mathcal U}_{ik} (\w_k )\qquad\mbox{for}\quad i=1,\ldots,p,
\]
with ${\mathcal U}_{ik}:L^{2}(\Omega,{\mathbb R}^{n})\rightarrow L^{2}(\Omega,{\mathbb R}^{n})$
 chosen so that, for $k=1,\ldots,i-1$,
\begin{gather}
\label{ewwik}
{\mathbb E}_{w_i w_k}=\mathbb O\quad\mbox{if}\quad i\neq k.
\end{gather}
We wish (\ref{ewwik}) is true for any $k$, i.e.
\begin{gather*}
{\mathbb E}_{w_i w_k}=E_{w_i w_k} - E[\w_i]E[\w_k^T] \\
\phantom{{\mathbb E}_{w_i w_k}}{} = E\Bigg[\Bigg(\uu_i - \sum_{l=1}^{i-1} {\mathcal U}_{il}
(\w_l)\Bigg)\w_k^T\Bigg] - E\Bigg[\Bigg(\uu_i - \sum_{l=1}^{i-1} {\mathcal U}_{il} (\w_l)\Bigg)\Bigg] E[\w_k^T]\\
\phantom{{\mathbb E}_{w_i w_k}}{} =  E_{u_i w_k} - U_{ik}E_{w_k w_k} - E[\uu_i] E[\w_k^T] + E[\w_k] E[\w_k^T]\\
\phantom{{\mathbb E}_{w_i w_k}}{} =  {\mathbb E}_{u_i w_k} -U_{ik}{\mathbb E}_{w_k w_k}=\mathbb O.
\end{gather*}
Thus, $U_{ik} = {\mathbb E}_{u_i w_k} {\mathbb E}_{w_k w_k}^{-1}$, and the statement (i) is true.

It is clear that vectors $\vv_1,\ldots,\vv_p$, def\/ined by (\ref{wq}), are orthogonal.
For ${\mathcal Q}_k$, def\/ined by (\ref{wq1}), we have $Q_k=({\mathcal E}_{w_k w_k}^{1/2})^{-1}$ and
\begin{gather*}
{\mathbb E}_{v_k v_k}
 = E\big[({\mathbb E}_{w_k w_k}^{1/2})^{-1} \w_k \w_k^T ({\mathbb E}_{w_k w_k}^{1/2})^{-1}\big]-
 E\big[({\mathbb E}_{w_k w_k}^{1/2})^{-1} \w_k]E[\w_k^T ({\mathbb E}_{w_k w_k}^{1/2})^{-1}\big]\\
\phantom{{\mathbb E}_{v_k v_k}}{}
=  ({\mathbb E}_{w_k w_k}^{1/2})^{-1} {\mathbb E}_{w_k w_k}({\mathbb E}_{w_k w_k}^{1/2})^{-1} = I.
\end{gather*}
Hence,  $\vv_1,\ldots,\vv_p$, def\/ined by (\ref{wq}), are orthonormal.

\medskip

\noindent
{\bf Proof of Lemma \ref{ort3}.} We wish that ${\mathbb E}_{v_i v_k} = {\mathbb O}$ for $i\neq k$.
If $Z_{ik}$ has been chosen so that this condition is true   for all $k=1,\ldots,i-1$ then
we have
\begin{gather}
\label{ev}
E\Bigg[\Bigg(\uu_i - \sum_{l=1}^{i-1} \z_{i l} (\vv_l)\Bigg)\vv_k^T\Bigg]
=  {\mathbb E}_{u_i v_k} - \sum_{l=1}^{i-1} Z_{il} {\mathbb E}_{v_l v_k}
 =  {\mathbb E}_{u_i v_k} - Z_{ik} {\mathbb E}_{v_k v_k}={\mathbb O}.
\end{gather}
Thus,
\begin{gather}
\label{auv}
Z_{i k}{\mathbb E}_{v_k v_k} = {\mathbb E}_{u_i v_k}.
\end{gather}
The necessary and suf\/f\/icient condition \cite{ben1} for the solution of the matrix equation  (\ref{auv}) is given by
\begin{gather}
\label{cond-eq}
{\mathbb E}_{u_i v_k}{\mathbb E}^\dag_{v_k v_k}{\mathbb E}_{v_k v_k}={\mathbb E}_{u_i v_k}.
\end{gather}
By Lemma \ref{phil}, (\ref{cond-eq}) is true. Then, on the basis of \cite{ben1}, the general  solution
 to (\ref{auv})  is given  by~(\ref{a}).

\subsection*{Acknowledgements}

 The f\/irst co-author  is grateful to Oliver Capp\'e {for useful discussions related to} the
 structure of the proposed transform.

\LastPageEnding

 \end{document}